\documentclass[10pt]{article}

\usepackage[english]{babel}
\usepackage{amsthm}
\usepackage{amssymb}
\usepackage{mathtools}

\newtheorem{theorem}{Theorem}[section]
\newtheorem{lemma}[theorem]{Lemma}
\newtheorem{proposition}[theorem]{Proposition}
\newtheorem{corollary}[theorem]{Corollary}
\newtheorem{remark}[theorem]{Remark}
\newtheorem{definition}[theorem]{Definition}
\makeatletter

\@addtoreset{equation}{section}
\usepackage{color}
\usepackage{graphicx}

\makeatother
\usepackage{newlfont}
\def \n {\noindent}
\begin{document}

\begin{center}
{\color{blue}{\bf {\large Review of some modified generalized Korteweg - de Vries - Kuramoto-Sivashinsky equations (mgKdV-KS) {\color{red} {\bf part I}}}}}\\
\end{center}

\begin{center}
{\bf Marie-Th\'er$\grave{e}$se Aimar$^{(1),(2)}$}\\

(1) CNRS, Centrale Marseille, I2M UMR 7373, Facult\'e des Sciences \quad \quad \quad \quad \\ Site St Charles Aix Marseille Universit\'e, \\3 place Victor Hugo, 13331 Marseille cedex 3, France\\

\n {\it (2) Le Prador,129, rue du commandant Rolland, 13008 Marseille, France}.\\
\n E.mail: {\color{blue}marietherese.aimar@orange.fr}\\
\quad\\

{\bf Abdelkader Intissar$^{(1)}$}\\
\n {\it (1) Le Prador, 129, rue du commandant Rolland, 13008 Marseille-France}.\\
\n E.mail: {\color{blue}abdelkader.intissar@gmail.com}\\
\end{center}

\begin{center}
{\color{red} {\bf Abstract}}
\end{center}
This paper reviews the results of existence and uniqueness of the solutions of these equations: the  Korteweg-de Vries equation, the Kuramoto-Sivashinsky equation, the generalized Korteweg-de Vries-Kuramoto-Sivashinski equation and the non homogeneous boundary value problem for KdV-KS equation in quarter plane. \\

\n {\bf Mathematics Subject Classification}. 35G25, 35K55.\\
\n {\bf Keywords}. Existence, Uniqueness, Stability, The modified generalized   Korteweg-de Vries-Kuramto-Sivashinsky equation, non linear initial value problems, Cauchy problem.\\
\begin{center}
{\color{blue} {\bf {\Large Content}}}
\end{center}
\n {\bf{\color{blue} $\S$\,\,0 }} {\color{black}Introduction} .......................................................................................... \hfill { } 2\\ 
\n {\bf{\color{blue} $\S$\,\,1 }} {\color{black}Elementary properties of Sobolev's spaces} .......................................... \hfill { } 5\\
\n {\bf{\color{blue} $\S$\,\,2 }} {\color{black}Some technicalities inequalities} ........................................................... \hfill { } 13\\ 
\n {\bf{\color{blue} $\S$\,\,3 }} {\color{black}Bonna-Smith results on Korteweg-de Vries (KdV) equation {\color{blue}(1975)}} ..... \hfill { }20\\
\n {\bf{\color{blue} $\S$\,\,4 }} {\color{black}Aimar's results on Kuramoto-Sivashinsky (K.S) equation {\color{blue}(1982)}} .......... \hfill { }26\\
\n {\bf{\color{blue} $\S$\,\,5 }} {\color{black} Presentation of KdV-KS equation {\color{blue}(1996)}} ............................................ \hfill { }32\\ 
\n {\bf{\color{blue} $\S$\,\,6}} {\color{black}\,  On travelling wave solution of the Kuramoto-Sivashinsky equation}  .... \hfill { }34\\
\n {\bf{\color{blue} $\S$\,\,7 }} {\color{black} Existence and uniqueness of the solution of KdV- KS equation on domain $[0, T ]\times [0, 2\pi], T > 0$}   .................................................................................... \hfill { } 35\\
\n {\bf{\color{blue} $\S$\,\,8 }} {\color{black}A study of the non homogeneous boundary value problem for KdV-KS equation in quarter plane {\color{blue} (2016)}}   .................................................................. \hfill { }47\\

 \n {\bf{\color{black}\Large{ 0 }}} {\bf{\color{red}\Large{Introduction}}}\\ 
 
\n In this paper , we review the results of the well-posedness of the classical solutions for the initial-boundary value problem for modified generalized Korteweg - de Vries-Kuramoto Shivashinsky equations (mgKdV-KS), under appropriate boundary conditions respectively:.We begin this work  by reviewing :\\

\n {\color{red}$\rhd$} In {\color{red}first section}, we present some properties of Sobolev spaces which are important there because they are the natural spaces to study  various partial differential equations  {\color{blue}[[35]]}, {\color{blue}[[29]]} . in particular to study modified generalized Korteweg - de Vries-Kuramoto Shivashinsky equations (mgKdV-KS).\\

\n {\color{red}$\rhd$} In {\color{red} second section}, we present  some technicalities inequalities which play a fundamental role for existence and uniqueness of the solutions of these equations  {\color{blue}[[9]]}.\\

\n {\color{red}$\rhd$} In {\color{red}Third section}, we consider the waves on the surface of an inviscid fluid in a flat channel. When one is interested in the propagation of one-directional irrotational small amplitude long waves, it is classical to model the waves by the well-known Korteweg - de Vries (KdV) equation.\\

\n {\color{red}$\bullet_{1}$}:  The expression of the Korteweg de Vries equation (KdV) equation {\color{blue}[[38]]},  {\color{blue} [[5]]}, {\color{blue} [[12]]},  {\color{blue} [[13]]},  {\color{blue} [[19]]} is given by:\\ 
\begin{equation}
\displaystyle{u_{t} + u_{x} + \epsilon uu_{x} + \mu u_{xxx} = 0}
 \end{equation}
\n or\\
 \begin{equation}
 \displaystyle{u_{t} +  uu_{x} +  u_{xxx} = f(x, t)}
 \end{equation}
 \n where $\displaystyle{u: \mathbb{R}^{2} \longrightarrow \mathbb{R}, (t, x) \longrightarrow u(t, x)}$ is a real function  and $u_{t}$ is its derivative on $t$, $u_{x}$ is its derivative on $x$ and $u_{xxx}$ is its theird derivative on $x$. $\epsilon$ and $\mu$ are the real parameters.\\
 
 \n {\color{red}$\rhd$} In {\color{red} forth section}, we consider the Kuramoto-Sivashinsky, which arises in the description of stability of flame fronts, reaction diffusion systems and many other physical settings  {\color{blue}[[24]]}, {\color{blue}[[33]]}, {\color{blue}[[21]]} {\color{blue}[[[2]] , [[3]]]}, {\color{blue}[[27]]} is one of the simplest nonlinear PDEs that exhibit spatiotemporally chaotic behavior. In the formulation adopted here, the time evolution of the flame front velocity $u = u(t, x)$ on a periodic domain $\displaystyle{u(t, x) = u(t, x + L) , L > 0} $.\\
 
 \n {\color{red}$\bullet_{2}$}:  The expression of the Kuramoto-Sivashinsky (KS) equation is given by:\\ 

\begin{equation}
\displaystyle{u_{t} + \frac{1}{2}(u_{x})^{2} + u_{xx} + u_{xxxx} = 0}
\end{equation}
\n or\:
\begin{equation}
\displaystyle{u_{t} + u u_{x} + u_{xx} + u_{xxxx} = 0}
\end{equation}
\n Here $t \geq 0$ is the time, and $x$ is the spatial coordinate. The subscripts $x$ and $t$ denote partial derivatives with respect to $x$ and $t$.\\

\n {\color{red}$\rhd$} In {\color{red}fifth section}, we consider the travelling wave solutions of the Kuramoto-Sivashinsky equation  of the form  $\tilde{u} = u(x - \omega t)$ and the steady-state solutions.
Substituting $\tilde{u} = u(\xi)$  where $\xi = x - \omega t$. Then the KS equation is transformed to \\
\begin{equation}
\displaystyle{- \omega u_{\xi} +  u_{\xi\xi\xi\xi}  + u_{\xi\xi}+ uu_{\xi} = 0 }
\end{equation}

\n This can be rewritten as a one-dimensional system by the change of variables $\displaystyle{x_{1} = \tilde{u}, x_{2} = \tilde{u}_{\xi}, x_{3} = \tilde{u}_{\xi\xi}}$ and $\displaystyle{ x_{4} = \tilde{u}_{\xi\xi\xi}}$\\

\n Then we deduce the Michelson system {\color{blue}[[27]] } \\

\begin{equation}
\left \{ \begin{array} {c} \displaystyle{\frac{dx_{1}}{d\xi} = x_{2}} \quad\quad\quad\quad\quad\quad\quad\quad\\
\quad\\
\displaystyle{\frac{dx_{2}}{d\xi} = x_{3}}  \quad\quad\quad\quad\quad\quad\quad\quad \\
\quad\\
\displaystyle{\frac{dx_{3}}{d\xi} = x_{4}}  \quad\quad\quad\quad\quad\quad\quad\quad\\
\quad\\
\displaystyle{\frac{dx_{4}}{d\xi} = \omega x_{2} - x_{1} x_{2}  - x_{3}\quad}\\
\end{array} \right.\\
\end{equation}

\n {\color{red}$\rhd$} In {\color{red}sixth section}, we consider the combined Korteweg-deVries equation with Kuramoto-Sivashinsky equation called generalized Korteweg-deVries-Kuramoto-sivashinsky equation (gKdV-KS) which arises in some interesting physical situations.\\ 

\n {\color{red}$\bullet_{3}$}:  The expression of the generalized Kuramoto-Sivashinsky (gKdV-KS) equation is given by:\\ 

\begin{equation}
\displaystyle{u_{t} + \frac{1}{2}(u_{x})^{2} +  u_{xxx}  + R u_{xx} + u_{xxxx} = 0}
\end{equation}
\n $R > 0$ is  Reynold number\\

\n Or the ``derivative'' generalized Korteweg-deVries-Kuramoto-sivashinsky equation\\
 \begin{equation}
\displaystyle{u_{t} + uu_{x} + \delta u_{xxx}  + \beta u_{xx} + \gamma u_{xxxx} = 0 }
\end{equation}
\n $ \delta, \beta, \gamma$ are constants.\\

\n \n  We recall that the term  $u_{xx}$ correspond to energy input at large scales (where the real positive number $R$ is called the ``anti-diffusion'' parameter), the term  $u_{xxxx}$ to dissipation at small scales and the term  $uu_{x}$ to nonlinear advection.\\

\n {\color{red}$\rhd$} In {\color{red}seventh section}, we consider  the local well-posedness and the global well-posedness of the non homogeneous initial
boundary value problem of KdV-KS equation in quarter plane styded by Jing Li,  Bing-Yu Zhang and Zhixiong Zhang  in {\color{blue} [[25]]}   .\\

\n In 2016, they was considered the following problem:\\
\begin{equation}
\left \{ \begin{array} {c} \displaystyle{u_t + u_{xxxx} + \delta u_{xxx} + u_{xx} + uu_{x} = 0, \quad (t, x) \in (0, T)\times \mathbb{R}^{+}}\\
\quad\\
u(0, x) = \phi(x) \quad\quad\quad \quad\quad\quad\quad \quad\quad \quad \quad\quad x\in \mathbb{R}^{+}\quad \quad\quad\quad \quad\\
\quad\\
u(t, 0) = h_{1}(t) \quad u_{x}(t, 0) = h_{2}(t) \quad \quad\quad \quad\quad t \in (0, T),\quad\quad \quad\quad\\
\end{array} \right. 
\end{equation}

\n where $T > 0$, $\displaystyle{\delta \in \mathbb{R}}$ and $\displaystyle{\phi \in H^{s}(\mathbb{R}^{+})}$.\\

\n They was conjectured that the wellposedness will fail when $s < -2.$.\\

\n Some their main results was based on $s$-compatibility conditions.\\

\n {\color{red}{\Large $\star$}} In another paper, we review some results concerned with the local well-posedness of the initial-value problems (IVP) for the Kawahara equation\\
\begin{equation}
 \left \{ \begin{array} {c} \displaystyle{u_{t} + \alpha u_{xxx} + \beta u_{xxxxx} +  uu_{x} = 0 \quad x, t \in \mathbb{R}}\\
  \quad\\
  u(0, x) = u_{0}(x) \quad \quad \quad\quad\quad \quad\quad \quad\quad \quad\quad \quad\\
  \end{array}  \right.
\end{equation}
and for the modified Kawahara equation
\begin{equation}
 \left \{ \begin{array} {c} \displaystyle{u_{t} + \alpha u_{xxx} + \beta u_{xxxxx} +  u^{2}u_{x} = 0 \quad x, t \in \mathbb{R}}\\
  \quad\\
  u(0, x) = u_{0}(x) \quad \quad \quad\quad\quad \quad\quad \quad\quad \quad\quad \quad\\
  \end{array}  \right.
\end{equation}
\n where $\alpha$ and $\beta$ are real constants and $\beta \neq  0$. These fifth-order KdV type equations arise in modeling gravity-capillary waves on a shallow layer and magneto-sound propagation in
plasmas (see e.g. {\color{blue}[[18]], [[20]]}).\\

\n The wellposedness issue on these fifth-order KdV type equations has previously been studied by several authors. In {\color{blue}[[30]]}, Ponce considered a general fifth-order KdV equation \\
\begin{equation}
 \left \{ \begin{array} {c} \displaystyle{u_{t} + c_{1} uu_{x} + c_{2} u_{xxx} + c_{3}u_{x}u_{xx} + c_{4}uu_{xxx} + c_{5} u_{xxxxx}  = 0 \quad x, t \in \mathbb{R}}\\
  \quad\\
  u(0, x) = u_{0}(x) \quad\quad\quad\quad \quad \quad\quad \quad \quad \quad \quad \quad \quad\quad\quad \quad\quad \quad\quad \quad\quad \quad\\
  \end{array}  \right.
\end{equation}

\n and established the global well-posedness of the corresponding IVP for any initial data in $H^{4}(\mathbb{R})$.\\

\section{{\color{red} {\bf Elementary properties of Sobolev's spaces}}}
  
\begin{definition} ({\color{red}Classical Lebesgue and Sobolev spaces})\\

\n Let $I$ be an interval of $\mathbb{R}$ and $p \in \mathbb{R}$, $1 \leq p \leq \infty$ we define the Sobolecv's space as following :\\

\n $\displaystyle{ W^{1,p}(I) = \{ u \in L^{p}(I) , \exists \,\, g \in L^{p}(I) ;  \int_{I}u(s) \varphi^{'}(s) ds = - \int_{I} g(s) \varphi(s)ds \quad \forall \,\, \varphi \in C_{c}^{1}(I) \}}$

\n where\\

\n $\displaystyle{L^{p}(I), 1 \leq p < \infty}$ : Space of the measurable functions $u$ on $I$ such that $\displaystyle{\int_{I}\vert u(s) \vert^{p}ds < \infty}$\\

\n $\displaystyle{\vert\vert u \vert\vert_{L^{p}(I)}, 1 \leq p < \infty}$: norm of $u$ in $L^{p}(I)$ given by $\displaystyle{\vert\vert u \vert\vert_{L^{p}(I)} = (\int_{I}\vert u(s) \vert^{p}ds)^{\frac{1}{p}}}$\\

\n  $C_{c}^{1}(I)$ is space of differentiable functions such that their derivative is continuous  with compact support.\\

\n We denote:\\

\n  $\displaystyle{H^{1}(I) = W^{1,{\color{red}2}}(I)}$.\\

\n  $\displaystyle{C^{n}(I)}$  Space of $n$-times differentiable functions defined on $I$ such that $k$-th derivative is continuous for $ k = 0, ....., n$\\

\n  $\displaystyle{L^{\infty}(I)}$: Space of measurable functions $u$ on $I$ such that are essentially bounded.\\

\n  $\displaystyle{\vert\vert u \vert\vert_{L^{\infty}(I)}}$: Norm of $u$ in the space$L^{\infty}(I)$ i.e, $\displaystyle{\vert\vert u \vert\vert_{L^{\infty}(I)} = sup \{\vert u(s) \vert ; s \in I\}}$.\\

\end{definition}

\begin{remark}

\n {\color{red} $\bullet_{1}$} $\displaystyle{ W^{1,p}(I) }$ can be defined by :\\
$$\displaystyle{ W^{1,p}(I)  = \{ u \in L^{p}(I) , u' \in L^{p}(I) \text{ (in weakly sense that correspond to} \,\, g \text{ of above definition}) \}}$$
\n {\color{red} $\bullet_{2}$} $\displaystyle{ W^{1,p}(I)}$ is equipped by  norm:\\
\begin{equation}
\displaystyle{\vert\vert u \vert\vert_{W^{1,p}(I)} = \vert\vert u \vert\vert_{L^{p}(I)} + \vert\vert u^{'} \vert\vert_{L^{p}(I)}}
\end{equation}

\n {\color{red} $\bullet_{3}$}  $\displaystyle{ H^{1}(I) }$ is equipped by scalar product :\\
\begin{equation}
\displaystyle{< u  , v >_{H^{1}(I)} = < u, v >_{L^{2}(I)} + < u^{'} , v^{'} >_{L^{2}(I)}}
\end{equation}
\n i.e\\
\begin{equation}
\displaystyle{ H^{1}(I) = \{u \in L^{2}(I) ; u' \in L^{2}(I)\};  \vert\vert u \vert\vert_{H^{1}(I)} = (\vert\vert u\vert\vert_{L^{2}}^{2} + \vert\vert u' \vert\vert_{L^{2}(I)}^{2})^{\frac{1}{2}}}
\end{equation}

\n {\color{red}$\bullet_{4}$} Let $k \in \mathbb{N}$, the Sobolev's space $H^{k}(I)$ is defined by:\\
\begin{equation}
\displaystyle{ H^{k}(I) = \{u \in L^{2}(I) ; u^{\alpha} \in L^{2}(I), \forall \,\, \alpha \leq k\}; \,\, \vert\vert u \vert\vert_{H^{k}(I)} = (\sum_{\alpha = 0}^{k}\vert\vert \partial^{\alpha}u\vert\vert_{L^{2}})^{\frac{1}{2}}}
\end{equation}
\n {\color{red}$\bullet_{5}$}  $\displaystyle{ W^{1,p}(I) }$ is Banach space for $1 \leq p \leq \infty$ and reflexive for  $1 < p < \infty$\\

\n {\color{red}$\bullet_{6}$} $\displaystyle{ H^{1}(I) }$ is Hilbert space with respect above scalar product.\\

\n {\color{red}$\bullet_{7}$} $\displaystyle{ H_{0}^{1}(I) }$ is the closure of $C_{c}^{1}(I)$ in $H^{1}(I) $.\\

\n {\color{red}$\bullet_{8}$} $\displaystyle{u \in L_{loc}^{1}(I)}$ $\iff$ $u \in L^{1}(K)$ for all compact $K \subset I$.

\end{remark}

\begin{theorem} ({\color{red}continuous representation})\\

\n{\color{red} (i)} Let $u \in  L_{loc}^{1}(I)$, if $\displaystyle{\int_{I}u\varphi = 0 \quad  \forall \varphi  \quad \in C_{c}(I)} $ then  $\displaystyle{u = 0 \,\, \text{ almost everywhere on} \,\,  I }$\\

\n{\color{red} (ii)} Let $u \in  L_{loc}^{1}(I)$ which satisfies $\displaystyle{\int_{I}u \varphi^{'} = 0 \quad \forall \,\, \varphi \in C_{c}^{1}(I)}$ then \\

\n $\displaystyle{\,\, \exists \,\, C \text{(a constant)} \,\,\text{such that }\,\, u = C\,\, \text{ almost everywhere}}$\\

\n{\color{red} (iii)} Let $u \in  L_{loc}^{1}(I)$, $x_{0} \in I$ and $\displaystyle{v(x) = \int_{x_{0}}^{x}g(s)ds \quad \forall x \in I}$ then \\

\n $v \in C(I)$ and $\displaystyle{\int_{I}v \varphi^{'} = -\int_{I}g \varphi}$\\

\n where $C(I)$ be the space of continuous functions on $I$.\\ 
 
\n{\color{red} (iv)} Let $u \in W^{1,p}(I)$ then  \\

\n  $\displaystyle{ \exists \,\, \tilde{u} \in C(I); u = \tilde{u}\,\, \text{almost everywhere on} \, I}$ and $\displaystyle{ \tilde{u}(x) -  \tilde{u}(y) = \int_{y}^{x}u'(s)ds \,\, \forall \,\, x, y \in I}$\\
\end{theorem}
\n {\color{red}{\bf Proof}}\\

\n{\color{red} (i)} Assume, at contrary, we suppose there exists a set $A \subset I$ such that $\mu(A) \neq 0$ and $u  > 0$ on $I$ (where $\mu$ is measure  of Lebesgue.\\

\n {\color{red} $\rhd_{1}$} If $I$ is bounded we put $a = inf I$ and $b = sup I$ and we consider the sequence $\displaystyle{([a + \frac{1}{n} , b - \frac{1}{n}])_{n}}$. $\displaystyle{\mu(A \cap [a + \frac{1}{n} , b - \frac{1}{n}]) \longrightarrow \mu(A) \neq 0}$ as $n \longrightarrow +\infty$, then there exists $n_{0} $ such that $\displaystyle{\mu(A \cap [a + \frac{1}{n_{0}} , b - \frac{1}{n_{0}}]) \neq 0}$.\\

\n {\color{red} $\rhd_{2}$} If $I$ is not semi-bounded  then we consider the sequence $\displaystyle{([-n , n])_{n}}$ and it follows  there exists $n_{0} $ such that $\displaystyle{\mu(A \cap [- n_{0} , n_{0}]) \neq 0}$.\\

\n {\color{red} $\rhd_{3}$} If $I$ has another form we can consider the following sequence $\displaystyle{([a - \frac{1}{n} , n])_{n}}$. So we can always find an interval $[c , d] \subset I$ such that $\displaystyle{\mu(A \cap ([c, d]) \neq 0}$.\\

\n We can take the following function $v$ such that \\

\n $\left \{ \begin{array} {c} v = 1 \,\, \text{ on}\,\, [c, d] \quad \quad  \quad \quad \\
\quad\\
v = 0 \,\, \text{ on}\,\, ^{C}[c - \frac{1}{n}, d + \frac{1}{n}]\\
\end{array} \right.$\\

\n From $v \in C_{c}(I)$ we deduce that $\displaystyle{\int_{I} uv = 0}$ and as $u > 0$ on $A$ it follows that $\displaystyle{\int_{I} uv \geq \int_{A\cap[c, d]} uv = \int_{A\cap[c, d]} u > 0}$ which is not possible . Then  $u = 0 $ almost everywhere on $ I$.\\

\n {\color{red} (ii)} Let $\psi \in  L_{loc}^{1}(I);  \int_{I} \psi = 1$.\\

\n $\displaystyle{\forall \,\, w \in C_{c}^{1}(I), \,\, \text{let}  h = w  - ( \int_{I} w)\psi \in C_{c}(I)}$ then $\exists \,\, \varphi \in C_{c}(I)$; $\displaystyle{\varphi^{'} = w  - ( \int_{I} w)\psi}$ it follows that\\

\n $\displaystyle{ \int_{I} u \varphi^{'} = 0}$\\

\n $\displaystyle{ \int_{I} u(x)( w(x)  - ( \int_{I} w(t)dt)\psi(x))dx   = 0 }$\\

\n $\displaystyle{ \int_{I} u(t)w(t)dt  - \int_{I} \int_{I} w(t)\psi(x)u(x)dtdx   = 0 }$\\

\n $\displaystyle{ \int_{I} w(t) (u(t)  - \int_{I}\psi(x)u(x)dx)dt   = 0 }$\\

\n As this last equality holds for all $w \in C_{c}(I)$ then by the property (i) we have $\displaystyle{u(t) -  \int_{I} \psi(x)u(x)dx   = 0 }$ almost everywhere on $I$\\

\n It follows that $u = C$  almost everywhere on $I$ where $\displaystyle{C = \int_{I} \psi(x)u(x)dx }$.\\

\n{\color{red} (iii)} Let $a = inf I$ and $b = sup I$ then we have\\

\n $\displaystyle{ \int_{I} v \varphi^{'} = \int_{a}^{b}(\int_{x_{0}}^{x}g(t)dt) \varphi^{'}(x)dx = \int_{a}^{x_{0}}\int_{x}^{x_{0}} g(t) \varphi^{'}(x)dtdx + \int_{x_{0}}^{b} \int_{x_{0}}^{x} g(t) \varphi^{'}(x)dt dx}$\\

\n By Fubini's theorem and as $\varphi(a) = \varphi(b) = 0$ then we obtain\\

\n $\displaystyle{ \int_{I} v \varphi^{'} = \int_{a}^{x_{0}} g(t) \int_{a}^{t}\varphi^{'}(x)dxdt + \int_{x_{0}}^{b}g(t)\int_{t}^{b} \varphi^{'}(x)dx dt = -  \int_{I}g\varphi}$\\

\n{\color{red} (iv)} Let a fixed $x_{0} \in I $ and $\displaystyle{ \tilde{u}(x) = \int_{x_{0}}^{x}u'(t)dt}$ then by (iii) we deduce that \\

$$\displaystyle{\int_{I}\tilde{u}\varphi^{'} = - \int_{I}u^{'} \varphi \quad \forall \,\, \varphi \in C_{c}^{1}(I)}$$

\n Using an integration by parts, we deduce that\\

$$\displaystyle{\int_{I}(u - \tilde{u})\varphi^{'} = 0 \quad \forall \,\, \varphi \in C_{c}^{1}(I)}$$
 
\n So by (iii) we deduce that $\displaystyle{u - \tilde{u} = C}$ almost  everywhere on $I$ and the function $\displaystyle{\tilde{u}(x) = u(x) + C}$ satisfies the request properties.\\

\begin{definition}
\quad\\

\n Let $\Omega$ an open domain of $\mathbb{R}^{n}$ and let $p \in \mathbb{R}$, $1 \leq p \leq \infty$. The Sobolev space $\displaystyle{ W^{1,p}(\Omega) }$ can be defined by :\\

\n {\color{red}$\bullet_{1}$}  $\displaystyle{W^{1,p}(\Omega) = \{u \in L^{p}(\Omega) , \exists \,\, g_{1}, g_{2}, ....., g_{n} \in  L^{p}(\Omega); \int_{\Omega}u\frac{\partial \varphi}{\partial x_{i}} = -\int_{\Omega}g_{i}\varphi \,\, \forall \,\, \varphi \in C_{c}^{1}(\Omega)\,\, \forall \,\, i \in [1, n]\}}$\\

\n or\\
$$\displaystyle{ W^{1,p}(\Omega)  = \{ u \in L^{p}(\Omega),\,\, \forall \,\, i \in [1, n]\,\, \frac{\partial u}{\partial x_{i}} \in L^{p}(\Omega) \text{ (in weakly sense that correspond to} \,\, g_{i}) \}}$$

\n We denote $\displaystyle{\nabla u =}$ $\left ( \begin{array} {c}  \frac{\partial u}{\partial x_{1}}\\
 \frac{\partial u}{\partial x_{2}}\\
.\\
.\\
.\\
\quad\\
 \frac{\partial u}{\partial x_{n}}\\
 \end{array} \right )$ = $\left ( \begin{array} {c}  g_{1}\\
g_{2}\\
.\\
.\\
.\\
\quad\\
 g_{n}\\
 \end{array} \right )$\\

\n {\color{red} $\bullet_{2}$} $\displaystyle{ W^{1,p}(\Omega)}$ is equipped by  norm:\\

\begin{equation}
\displaystyle{\vert\vert u \vert\vert_{W^{1,p}(\Omega)} = \vert\vert u \vert\vert_{L^{p}(\Omega)} + \vert\vert \nabla u \vert\vert_{L^{p}(\Omega)}}
\end{equation}

\n {\color{red} $\bullet_{3}$} We denote $\displaystyle{ H^{1}(\Omega) = W^{1,2}(\Omega)}$ which is equipped by scalar product :\\

\begin{equation}
\displaystyle{< u  , v >_{H^{1}(\Omega)} = < u, v >_{L^{2}(\Omega)} + \sum_{i=1}^{n}<  \frac{\partial u}{\partial x_{i}} , \frac{\partial v}{\partial x_{i}} >_{L^{2}(\Omega)}}
\end{equation}

\n associated to the following norm:\\

\begin{equation}
\displaystyle{ \vert\vert u \vert\vert_{H^{1}(\Omega)} = (\vert\vert u\vert\vert_{L^{2}(\Omega)}^{2} + \vert\vert \nabla u \vert\vert_{L^{2}(\Omega)}^{2})^{\frac{1}{2}}}
\end{equation}

\n {\color{red} $\bullet_{4}$} $\displaystyle{ H_{0}^{1}(\Omega)}$ is defined as the closure of $C_{c}^{1}(\Omega)$ in $H^{1}(\Omega)$.
\end{definition}
\begin{remark}
\n There is no continuous representation like that one dimension of above theorem\\
\end{remark}

\begin{definition}
\n Let $m \geq 2$ and  $1 \leq p \leq \infty$, by recurrence we define:\\\begin{equation}
\displaystyle{W^{m,p}(\Omega) = \{u \in W^{m-1,p}(\Omega), \frac{\partial u}{\partial x_{i}} \in W^{m-1,p}(\Omega) \,\, \forall  \,\, i \in [1, n]\}}
\end{equation}
\n and\\
\begin{equation}
\displaystyle{H^{m}(\Omega) = W^{m,2}(\Omega)} 
\end{equation}
\end{definition} 
\begin{remark}
\quad\\

\n {\color{red}$\rhd_{1}$} The Sobolev spaces {\color{blue}[[35]]} were introduced mostly for the use of the theory of partial differential equations. Differential operators are often closable in such spaces. Of course, Sobolev spaces being examples of Banach or, sometimes, Hilbert spaces are interesting object for themselves. But their importance is connected with the fact that the theory of partial differential equations can be, and even most easily, developed just in such a spaces. The reason is because partial differential operators are very well situated in Sobolev spaces.. The spaces of continuous (or of class $C^{k}$) functions is not very suitable for the studies of partial differential equations.\\

\n Why spaces of (only) continuous functions are not very suitable? The answer is connected with the following observation  {\color{blue}[[29]]} \\

\n  Namely, for every $k \in \mathbb{N}$ the Laplace operator
\begin{equation}
\displaystyle{   \Delta ; C^{k+2}(\Omega) \longrightarrow  C^{k}(\Omega) ; u = 0 \,\, \text{on} \,\, \partial \Omega}
\end{equation}
\n is continuous, but its image is not closed in $C^{k}(\Omega)$. In particular, for continuous right hand side $f \in C^{k}(\Omega)$ the solution of the equation:\\
\begin{equation}
\displaystyle{   \Delta u = f, \quad u = 0 \,\, \text{on} \,\, \partial \Omega}
\end{equation}
in general must not be a  $C^{2}(\Omega)$ function. Similar situation we face for other elliptic operators.\\

\n {\color{red}$\rhd_{2}$} Consider now a partial differential operator with constant coefficients $a_{\alpha}$:\\

\begin{equation}
\displaystyle{ P(D) = \sum_{\vert \alpha\vert \leq m}a_{\alpha}D^{\alpha}}
\end{equation}
\n It can be considered as an operator from $C^{\infty}(\Omega)$ into itself, or as an operator from $C^{k}(\Omega)$ into $C^{k - m}(\Omega)$, for any $m  \leq k$. But we prefer to define it as:\\

\begin{equation}
\displaystyle{ P(Du)(x) = \sum_{\vert \alpha\vert \leq m}a_{\alpha}D^{\alpha}u(x)}
\end{equation}
\n on the linear subspace $\mathfrak{C}^{m}(\Omega)$ consisting of all the functions in $C^{m}(\Omega)$ that have a finite norm:\\
\begin{equation}
\displaystyle{ \vert\vert u \vert\vert_{W^{m, p}(\Omega)} = (\sum_{\vert \alpha \vert \leq m} \int_{\Omega} \vert D^{\alpha}u \vert^{p})^\frac{1}{p}}
\end{equation}

\n It is often desirable to extend $P(D)$ to a closed linear operator in $L^{p}(\Omega),(1\leq p < \infty)$. It is well known that this is possible:\\
\end{remark}
\begin{theorem}
\quad\\

\n The operator $P(D)$ from $ L^{p}(\Omega)\,\,\, p \in [1; \infty) $ into itself with a domain $\mathfrak{C}^{m}(\Omega)$ has a closure.\\
\end{theorem}
\n {\color{red}{\bf Proof}}\\

\n The following is a well known equivalent condition for  an operator $T$ to have a closure:\\

\n {\color{red}$\rhd$}  Let $T$ be a linear operator from a linear subspace $D(T)$ of a Banach space $\mathbb{X}$ into a Banach space $\mathbb{Y}$ . $T$ has a closure $\tilde{T}$ if and only if the following condition is satisfied:\\
 \begin{equation}
 \displaystyle{u_{n} \in D(T) ,  u_{n} \longrightarrow 0, Tu_{n}  \longrightarrow v \,\, \text{ imply}\,\, v = 0}
 \end{equation}
 
  \n {\color{red}$\rhd$} Now, in view of the necessary condition above it suffices to show that if $\displaystyle{u_{n} \in C^{m}(\Omega); \vert\vert u_{n} \vert\vert_{W^{0,p}}  \longrightarrow 0, \vert\vert P(D)u_{n} - v\vert\vert_{W^{0,p}}  \longrightarrow 0}$ then $v = 0$.\\
 
 \n Let $\varphi \in C_{0}^{\infty}(\Omega)$. Integration by parts gives:\\
 
 \begin{equation}
 \displaystyle{\int_{\Omega}P(D)u_{n}\varphi dx = {\int_{\Omega}u_{n} P(-D)}\varphi dx}
 \end{equation}
 
 \n As $n \longrightarrow \infty$, the integrals on the right converge to zero whereas the integrals on the left converge to $\displaystyle{\int_{\Omega} v\varphi dx}$. Therefore,\\
 
 \begin{equation}
 \displaystyle{\int_{\Omega} v\varphi dx = 0 \,\, \text{for all} \,\, \varphi \in C_{0}^{\infty}(\Omega)}
 \end{equation} 
 
 \n Thus, by density, we conclude that $v = 0$ and this result show that the Sobolev spaces are natural for the studies of differential operators.\\
 
\n  For more informations concerning dense subsets of Sobolev spaces  and the classes of domains, described in terms of the ``smoothness'' of their boundary $\partial \Omega$ \\

\n considered in this Sobolev theory consult ({\color{blue}[[1]]}, Chapter III).\\ In particular the three classes of domains which are most often considered:\\

\n $\bullet$ Domains $\Omega \subset \mathbb{R}^{n}$ having the {\color{red}cone property}.\\

\n $\bullet$ Domains having the {\color{blue}local Lipschitz property}.\\

\n $\bullet$ Domains having the {\color{green}$C^{m}$-regularity property}.\\

\n {\color{red}{\bf The fractional Sobolev space $\displaystyle{W^{s, p}}$}}: Let $\Omega$ be a general open set in $\mathbb{R}^{n}$. For any real $s > 0$ and for any $p \in  [1,  +\infty[$, we want to define the {\color{red}fractional Sobolev spaces}  $\displaystyle{W^{s, p}(\Omega)}$. In the literature, fractional Sobolev-type spaces are also called Aronszajn, Gagliardo or Slobodeckij spaces, by the name of
the ones who introduced them, almost simultaneously (see {\color{blue} [[4]], [[11]] , [[34]]}).\\

\n We start by fixing the fractional exponent $s$ in $]0, 1[$. For any $p \in [1, + \infty[$, we define $\displaystyle{W^{s, p}(\Omega)}$ as follows \\

\begin{equation}
\displaystyle{W^{s, p}(\Omega) = \{u \in L^{p}(\Omega) ; \frac{\vert u(x) - u(y)\vert^{p}}{\vert x - y\vert^{\frac{n}{p} + s}} dxdy \in L^{p}(\Omega \times \Omega)}
\end{equation}

\n i.e., an intermediary Banach space between $L^{p}(\Omega)$ and $W^{1,p}(\Omega)$ endowed with the natural norm :\\

\begin{equation}
\displaystyle{\vert\vert u \vert\vert_{W^{s, p}(\Omega)} := (\int_{\Omega} \vert u \vert ^{p} +  \int_{\Omega}\int_{\Omega} \frac{\vert u(x) - u(y)\vert^{p}}{\vert x - y\vert^{n + ps}}dx dy)^{\frac{1}{p}}}
\end{equation}
\n where the term\\
$$\displaystyle{[ u ]_{W^{s, p}(\Omega)} := (\int_{\Omega}\int_{\Omega} \frac{\vert u(x) - u(y)\vert^{p}}{\vert x - y\vert^{n + ps}}dx dy)^{\frac{1}{p}}}$$

\n is the so-called {\color{red}Gagliardo semi-norm} of $u$.\\

\begin{proposition}
.\n  Let $p \in [1,+ \infty[$ and $0 < s \leq {\color{red}s^{'}} < 1$. Let $\Omega$ be an open set in $\mathbb{R}^{n}$ and $u : \Omega \longrightarrow \mathbb{R}$ be a measurable function. Then \\
\begin{equation}
\displaystyle{ \vert\vert u \vert\vert_{W^{s, p}(\Omega)} \leq \vert\vert u \vert\vert_{W^{ {\color{red}s^{'}} , p}(\Omega)} }
\end{equation}
\n for some suitable positive constant $C = C(n, p, s) \geq 1$. In particular \\
$$ \displaystyle{W^{ {\color{red}s^{'}} , p}(\Omega) \subset W^{s, p}(\Omega)}$$

\end{proposition}
 \n {\color{red}{\bf Proof}}\\
 
\n Since $n + sp > n$ then the kernel $\displaystyle{\frac{1}{\vert z \vert^{n + ps}}}$  is integrable. It follows that\\

\n $\displaystyle{\int_{\Omega} \, \int_{\Omega \cap \{\vert x - y\vert \geq 1\}} \frac{\vert u(x) \vert^{p}}{\vert x - y\vert^{n + sp}}dxdy \leq \int_{\Omega}(\int_{\vert z \vert \geq 1} \frac{1}{\vert z \vert
^{n + ps}}dz) \vert u(x) \vert^{p}dx}$ \\

\n $\displaystyle{\leq C(n,p,s) \vert\vert u \vert\vert_{L^{p}(\Omega)}^{p}}$.\\

\n Taking into account the above estimate, it follows \\

\n $\displaystyle{\int_{\Omega} \, \int_{\Omega \cap \{\vert x - y\vert \geq 1\}} \frac{\vert u(x) - u(y)\vert^{p}}{\vert x - y\vert^{n + sp}}dxdy \leq 2^{p-1}\int_{\Omega} \, \int_{\Omega \cap \{\vert x - y\vert \geq 1\}} \frac{\vert u(x) \vert^{p} + \vert u(y) \vert^{p}}{\vert x - y\vert^{n + sp}}dxdy }$
\begin{equation}
\displaystyle{ \leq 2^{p}C(n, p, s) \vert\vert u \vert\vert_{L^{p}(\Omega)}^{p}}
\end{equation}
\n On the other hand,\\

\begin{equation}
\displaystyle{\int_{\Omega} \, \int_{\Omega \cap \{\vert x - y\vert < 1\}} \frac{\vert u(x) - u(y)\vert^{p}}{\vert x - y\vert^{n + sp}}dxdy \leq \int_{\Omega} \, \int_{\Omega \cap \{\vert x - y\vert < 1\}} \frac{\vert u(x) - u(y)\vert^{p}}{\vert x - y\vert^{n + {\color{red}s^{'}}p}}dxdy}
\end{equation}

\n Thus, combining the two last inequalities, we get \\
\n $\displaystyle{\int_{\Omega} \, \int_{\Omega} \frac{\vert u(x) - u(y)\vert^{p}}{\vert x - y\vert^{n + sp}}dxdy \leq 2^{p}C(n, p, s) \vert\vert u \vert\vert_{L^{p}(\Omega)}^{p} + \int_{\Omega} \, \int_{\Omega} \frac{\vert u(x) - u(y)\vert^{p}}{\vert x - y\vert^{n +  {\color{red}s^{'}}p}}dxdy }$\\

\n and so\\

\n  $\displaystyle{\vert\vert u \vert\vert_{W^{s, p}(\Omega)}  \leq (C(n, p, s)  + 1) \vert\vert u \vert\vert_{L^{p}(\Omega)}^{p} + \int_{\Omega} \, \int_{\Omega} \frac{\vert u(x) - u(y)\vert^{p}}{\vert x - y\vert^{n +  {\color{red}s^{'}}p}}dxdy }$\\

\n i.e.\\

\n  $\displaystyle{\vert\vert u \vert\vert_{W^{s, p}(\Omega)}  \leq \tilde{C}(n, p, s) \vert\vert u \vert\vert_{W^{{\color{red}s^{'}}, p}(\Omega)}}$\\

\n which gives the desired estimate.\\

\begin{remark}
It is also possible to define spaces $W^{s,p}(\Omega)$ for $s \in \mathbb{R}-\mathbb{N}, s > 1$. This requires appearing the derivatives of integer orders less than s in the
Gagliardo semi-norm, as presented in {\color{blue}[[7, chapter 4.6]] or [[8]]}.\\
\end{remark}

\n In the following sections, we frequently use following inequalities.\\
 \section{{\color{red} { \bf Some technicalities inequalities}}}
 \begin{lemma}
 $\displaystyle{ \forall x, y \geq 0 , \forall \,\, 0 \leq \alpha \leq 1}$ then we have\\
 \begin{equation}
 \displaystyle{x^{\alpha}y^{1-\alpha} \leq \alpha x + (1 - \alpha)y}
 \end{equation}
 \end{lemma}
 
 \n {\color{red}{\bf Proof}}\\
 
 \n If $x = y = 0$, then the inequality holds. So we assume that $(x, y) \neq (0, 0)$  and as the inequality is symmetric with respect  $\alpha \longleftrightarrow (1 - \alpha)$ and $x \longleftrightarrow y$, we can assume that $y \neq 0$.\\
 
 \n Now since the set of couples $(x, y) \in \mathbb{R}_{+}\times \mathbb{R}_{+}^{*}$ is the same as the set of couples ($xt, y)$ with $y \neq 0$, we are reduced to establishing:\\
 \begin{equation}
 \displaystyle{(xy)^{\alpha}y^{1-\alpha} \leq \alpha xy + (1 - \alpha)y}
 \end{equation}
\n which, after division by y, becomes:\\
 \begin{equation}
 \displaystyle{x^{\alpha} \leq^{{\color{red}?}} \alpha x + (1 - \alpha) }
 \end{equation}
\n Since the inequality to be established is obvious for $\alpha = 0, 1$, we can assume that $ 0 < \alpha < 1$. We are therefore to consider here the function of a single variable $x \in \mathbb{R}_{+}$ :\\
 \begin{equation}
 \displaystyle{f(x) = x^{\alpha} - \alpha x +  \alpha - 1 }
 \end{equation}
So its derivative is:\\
 \begin{equation}
 \displaystyle{f'(x) = \alpha(x^{\alpha-1} - 1)} 
 \end{equation}
 \n Now we have $f(0) = \alpha -1 < 0$, $f'(x) \geq 0$ for $0 \leq x \leq 1$ and $f'(x) \leq 0$ for $1 \leq x \leq \infty$ which forces $f$ to reach its maximum at the point $x = 1$ then  as $f(1) = 0$ and therefore $f$ always takes values $\leq 0$, which establishes the desired inequality.\\
 
 \n From this lemma we deduce\\
 \begin{theorem} {\color{red} Holder inequality})\\
 
 \n Let $1 < p  < \infty$ and $q$ its  conjugate exponent i.e $\displaystyle{\frac{1}{p} + \frac{1}{q} = 1}$.\\
 
 \n Let $ u \in L^{p}(\mathbb{R}^{n}) $ and $ v \in L^{p}(\mathbb{R}^{n}) $ then we have\\
 \begin{equation}
 \displaystyle{ \vert\vert uv \vert\vert_{ L^{1}(\mathbb{R}^{n})} \leq \vert\vert u \vert\vert_{ L^{p}(\mathbb{R}^{n})} \vert\vert v \vert\vert_{ L^{q}(\mathbb{R}^{n})} \,\,( uv \in L^{1}(\mathbb{R}^{n}))}
 \end{equation}
 
 \n In particular\\
  \begin{equation}
 \displaystyle{\vert\vert uv \vert\vert_{ L^{1}(\mathbb{R}^{n})} \leq \vert\vert u \vert\vert_{ L^{2}(\mathbb{R}^{n})} \vert\vert v \vert\vert_{ L^{2}(\mathbb{R}^{n})}\,\, ( if\,\, p = q =2 )}
 \end{equation}
 and\\
  \begin{equation}
  \displaystyle{\vert\vert uv \vert\vert_{ L^{1}(\mathbb{R}^{n})} \leq \vert\vert u \vert\vert_{ L^{\infty}(\mathbb{R}^{n})} \vert\vert v \vert\vert_{ L^{1}(\mathbb{R}^{n})} \,\, (if \,\,p =\infty \,\, \text{therefore}\,\,  q = 1)}
 \end{equation}
 \end{theorem}
 \n{\color{red}{\bf Proof}}\\
 
 \n We can assume that $\displaystyle{\vert\vert u \vert\vert_{ L^{p}(\mathbb{R}^{n})} \neq 0}$ and $\displaystyle{\vert\vert v \vert\vert_{ L^{q}(\mathbb{R}^{n})} \neq 0}$ and divide $u$ and $v$ by their norms:\\
 
$$\displaystyle{u \longrightarrow \frac{u}{\vert\vert u \vert\vert_{ L^{p}(\mathbb{R}^{n})}}} \,\, \text{and}\,\,\displaystyle{v \longrightarrow \frac{v}{\vert\vert v \vert\vert_{ L^{p}(\mathbb{R}^{n})}}}$$

 \n in order to bring us back, in the inequality to be established, to the case where $u$ and $v$ are both of unit norm:\\
 
  $$\displaystyle{\vert\vert u \vert\vert_{ L^{p}(\mathbb{R}^{n})} = 1}\,\, \text{ and}\,\, \displaystyle{\vert\vert v \vert\vert_{ L^{q}(\mathbb{R}^{n})} = 1}$$
 
 \n For $t \in \mathbb{R}^{n}$ applying above lemmas to numbers $x = \vert u(t) \vert^{p}$ and  $y = \vert v(t) \vert^{q}$ with $\alpha = \frac{1}{p}$ and $1 - \alpha = \frac{1}{q}$ to get:\\
 
  $$\displaystyle{\vert u(t)\vert \vert v(t)\vert \leq \frac{1}{p}\vert u(t) \vert^{p} + \frac{1}{q}\vert v(t) \vert^{q}}$$
  
  \n Or \\
  
  \n {\color{red}$\bullet_{1}$} {\color{blue}Young's inequality :} $\displaystyle{ \forall a, b \geq 0 \,\,, \forall \epsilon > 0, \quad 1 < p < \infty , \frac{1}{p} + \frac{1}{q} = 1}$\\
  \begin{equation}
 \displaystyle{ab \leq \frac{a^{p}}{p} + \frac{b^{q}}{q}}
\end{equation}

  \n {\color{red}$\bullet_{2}$} (Young's inequality with {\color{red}$\epsilon$} {\color{blue} [[9]]}). \\

\n Let $\displaystyle{a > 0, b > 0, 1< p, q <\infty; \frac{1}{p} + \frac{1}{q} = 1}$, Then\\
\begin{equation}
\displaystyle{ ab \leq \epsilon a^{p}+ C_{\epsilon}b^{q}}, \quad \text{for} \quad \displaystyle{C_{\epsilon} = (\epsilon p)^{-q/p} q^{-1}}
\end{equation}

 \n So\\
 
  $$\displaystyle{\int_{\mathbb{R}^{n}}\vert u(t)\vert \vert v(t)\vert dt \leq \frac{1}{p}\int_{\mathbb{R}^{n}}\vert u(t) \vert^{p}dt + \frac{1}{q}\int_{\mathbb{R}^{n}}\vert v(t) \vert^{q}dt}$$
  
 \n 0r\\
 
   $$\displaystyle{\int_{\mathbb{R}^{n}}\vert u(t)\vert \vert v(t)\vert dt \leq \frac{1}{p}\vert\vert u \vert\vert_{ L^{p}(\mathbb{R}^{n})}^{p} + \frac{1}{q}\vert\vert v \vert\vert_{ L^{q}(\mathbb{R}^{n})}^{q} = \frac{1}{p} + \frac{1}{q} = 1}$$
   
\n i.e\\   
   
$$\displaystyle{\vert\vert uv \vert\vert_{ L^{1}(\mathbb{R}^{n})} \leq 1 = 1^{2} = \vert\vert u \vert\vert_{ L^{p}(\mathbb{R}^{n})}.\vert\vert v \vert\vert_{ L^{q}(\mathbb{R}^{n})}}$$

\begin{corollary} ({\color{red} interpolation inequalities for $L^{p}(\Omega)$ and for $ W^{m,p}(\Omega)$})\\

\n As a consequence of the Holder inequality the following inequalities are satisfied :\\
{\color{red}$\rhd_{1}$} \begin{equation}
\displaystyle{\vert\vert u\vert\vert_{ L^{q}(\Omega)} \leq \vert\vert u\vert\vert_{ L^{o}(\Omega)}^{\beta} \vert\vert u\vert\vert_{ L^{r}(\Omega)}^{1 - \beta}} 
\end{equation}
\n valid for $\displaystyle{u \in  L^{q}(\Omega)}$ with $\displaystyle{ p \leq q \leq r }$ and $\displaystyle{\frac{1}{q} = \frac{\beta}{p} + \frac{1- \beta}{r}}$

\n {\color{red}$\rhd_{2}$} \begin{equation}
\displaystyle{\vert\vert u\vert\vert_{j,p}  \leq K\epsilon \vert\vert u\vert\vert_{m,p}  + K\epsilon^{- \frac{j}{m-j}}\vert\vert u\vert\vert_{ 0,p} }
\end{equation}
for any $\displaystyle{u \in W^{m,p}(\Omega)}$ where $\Omega$ is a suitably regular domain.\\

\n Here $\displaystyle{0 \leq j \leq m-1}$ and \\

$$ \displaystyle{\vert\vert u\vert\vert_{j,p}  = \{\sum_{\vert \alpha \vert = j}\int_{\Omega}\vert D^{\alpha}u \vert^{p} \}^{\frac{1}{p}}}$$

\n where $\displaystyle{\alpha = ( \alpha_{1}, \alpha_{2}, ...., \alpha_{n}) \in \mathbb{N}^{n}, \vert \alpha \vert = \alpha_{1} + \alpha_{2} + ....+ \alpha_{n}}$ and  $\displaystyle{D^{\alpha}u = D_{1}^{\alpha_{1}}D_{2}^{\alpha_{2}} ....D_{n}^{\alpha_{n}}u}$
\end{corollary}

\n {\color{red}$\bullet$}  Can relating the integral of a function to the integral of its gradient ?\\

\begin{lemma} ({\color{red}The Dirichlet Poincar\'e Inequality for the case $n= 1$})\\

\n If $\displaystyle{u : ]-r , r[ \longrightarrow \mathbb{R}, r > 0}$ is a $C^{1}$ function satisfying $u(-r) = u(r) = 0$.then\\

\begin{equation}
\displaystyle{\int_{-r}^{r}u^{2}dx \leq 4r^{2}\int_{-r}^{r}(u')^{2}dx}
\end{equation}
\end{lemma}

\n {\color{red}{\bf Proof}}\\

\n By the Fundamental Theorem of Calculus we have $\displaystyle{u(s) = \int_{-r}^{s}u'(x)dx}$ Therefore $\displaystyle{\vert u(s) \vert \leq \int_{-r}^{s}\vert u'(x) \vert dx}$.\\

\n Applying the CauchySchwarz inequality to $\displaystyle{ \int_{-r}^{s}u'(x)dx}$ to get:  \\

\n $\displaystyle{ \int_{-r}^{s}u'(x)dx \leq ( \int_{-r}^{s}1^{2}dx)^{1/2}( \int_{-r}^{s}(u'(x))^{2}dx)^{1/2}}$\\

 \n $\displaystyle{\leq (s + r)^{\frac{1}{2}} \int_{-r}^{s}(u'(x))^{2}dx)^{1/2} }$\\
 
 \n $\displaystyle{ \leq (2r)^{1/2} \int_{-r}^{r}(u'(x))^{2}dx)^{1/2} }$.

\n It follows that:\\

\n $\displaystyle{\vert u(s) \vert \leq  (2r)^{1/2} \int_{-r}^{r}(u'(x))^{2}dx)^{1/2} }$.

\n Squaring both sides gives\\

\n $\displaystyle{\vert u(s) \vert^{2} \leq  2r \int_{-r}^{r}(u'(x))^{2}dx }$.\\

\n and finally we integrate over $[-r, r]$ to give\\

\n $\displaystyle{\int_{-r}^{r}u^{2}dx \leq 4r^{2}\int_{-r}^{r}(u')^{2}dx}$\\ 

\n It is not difficult to extend these proofs to higher dimensional cubes \\$B_{r} \subset \mathbb{R}^{n}, r > 0$.\\

\begin{corollary}
\n If $\displaystyle{u : B_{r}  \longrightarrow  \mathbb{R}}$ is a $C^{1}$ function with $u = 0$ on $\partial B_{r}$ then there exists a constant $C(n) > 0$ such that\\

\begin{equation}
\displaystyle{\int_{B_{r}} u^{2}dx \leq C(n)r^{2}\int_{B_{r}} \vert \nabla u\vert^{2}dx}
\end{equation}
\end{corollary}

\begin{lemma} ({\color{red}The Nueman-Poincar\'e Inequality for the case $n= 1$})\\

\n If $\displaystyle{u : ]-r , r[ \longrightarrow \mathbb{R}, r > 0}$ is a $C^{1}$ function  and we define $\displaystyle{A = \frac{1}{2r}\int_{-r}^{r}u dx}$ then\\

\begin{equation}
\displaystyle{\int_{-r}^{r}(u - A)^{2}dx \leq 4r^{2}\int_{-r}^{r}(u')^{2}dx}
\end{equation}
\end{lemma}

\n {\color{red}{\bf Proof}}\\

\n Take a differentiable function $u$ with $\displaystyle{A = \frac{1}{2r}\int_{-r}^{r}u dx}$. Note by the intermediate value theorem that there is a point $c$ in $[-r, r] $ with $u(c) = A$. We have\\

$$\displaystyle{u(s) = u(c) + \int_{c}^{s}u'(x)dx}$$

\n From this we get\\

$$\displaystyle{\vert u(s) - A \vert \leq  \int_{c}^{s}\vert u'(x)\vert dx \leq \int_{-r}^{r}\vert u'(x)\vert dx}$$

\n Apply Cauchy-Schwarz again to give\/

$$\displaystyle{\vert u(s) - A \vert \leq  2r (\int_{-r}^{r}(u'(x))^{2} dx)^{\frac{1}{2}}}$$

\n Squaring and integrating then gives :\\

$$\displaystyle{\int_{-r}^{r}\vert u(s) - A \vert^{2}ds \leq  4r^{2} \int_{-r}^{r}(u'(s))^{2} ds}$$

\n It is not difficult to extend these proofs to higher dimensional cubes\\ $B_{r} \subset \mathbb{R}^{n}, r > 0$.\\

\begin{corollary}
\n If $\displaystyle{u : B_{r}  \longrightarrow  \mathbb{R}}$ is a $C^{1}$ function and we define \\

$$\displaystyle{ A = \frac{1}{Vol(B_{r})}\int_{B_{r}} udx }$$ then there exists a constant $C(n) > 0$ such that\\

\begin{equation}
\displaystyle{\int_{B_{r}} (u - A)^{2}dx \leq C(n)r^{2}\int_{B_{r}} \vert \nabla u\vert^{2}dx}
\end{equation}
\end{corollary}

\begin{remark} ({\color{red}Poincar\'e inequality on $\mathbb{R}$})\\

\n The Poincar\'e's inequality can not hold on the unbounded domain  $\mathbb{R}$.\\ 
Indeed, consider the sequence of smooth functions \\

$\phi_{k}(x) = $ $\left \{ \begin{array} {c} 0 ; \quad \vert x \vert > k + \frac{1}{10}\\
\quad\\
- sign(x);\,\, \vert x \vert \in (k , k +1)\\
\quad\\
0 ; \quad \vert x \vert < k - \frac{1}{10}\\
\end{array} \right. $

\n where the function is smooth and monotone. Then for all $1\leq p < \infty $ we have $\vert\vert \phi_{k} \vert\vert_{L^{p}(\mathbb{R})} \equiv 2$ for all $k$ while the smooth functions  $\displaystyle{\psi(x) := \int_{-\infty}^{x}\phi_{k}(s)ds}$ satisfy $\vert\vert \psi_{k} \vert\vert_{L^{p}(\mathbb{R})} \longrightarrow \infty$   as $k \longrightarrow + \infty$. Thus, it is not possible to find a constant $C > 0$ such
that $\vert\vert \psi_{k} \vert\vert_{L^{p}(\mathbb{R})} \leq  C\vert\vert \psi_{k}^{'} \vert\vert_{L^{p}(\mathbb{R})}$ for all $k$, and hence the Poincar\'e inequality must fail in $\mathbb{R}$.\\
\end{remark}

\begin{lemma} ({\color{red}Poincar\'e inequality for the case $n \geq 2$})\\

\n Let  $\Omega \subset \mathbb{R}^{n}$ be a domain with continuous boundary. Let $1 \leq p < \infty$. Then for any $u \in W^{1, p} (\Omega)$ :\\
\begin{equation}
\displaystyle{\int_{\Omega} \vert u - \frac{1}{\vert \Omega \vert}\int_{\Omega} u(y)dy\vert^{p}dx \leq c \int_{\Omega} \sum_{i=1}^{n}\vert \frac{\partial u}{\partial x_{i}} \vert^{p}dx}
\end{equation}
\end{lemma}

\n {\color{red}{\bf Proof}}\\

\n The proof is equivalent with showing that:\\

$$\displaystyle{\int_{\Omega} \vert u \vert^{p}dx \leq c \int_{\Omega} \sum_{i=1}^{n}\vert \frac{\partial u}{\partial x_{i}} \vert^{p}dx}$$

\n for the $W^{1,p}(\Omega)$ functions having zero main value.\\

\n We will used  the following result:\\

\begin{theorem} (consult {\color{blue} [[1]]})\\

\n  The embedding $\displaystyle{W^{1,p}(\Omega) \subset L^{p}(\Omega)}$ is compact.\\
\end{theorem}

\n Now, assume, at contrary, that there exists a sequence $u_{k} \in W^{1,p}(\Omega), k \in \mathbb{N} $ such that $\displaystyle{\vert\vert u_{k} \vert\vert_{ W^{1,p}(\Omega)} = 1}$ and\\

\begin{equation}
\displaystyle{1 = \vert\vert u_{k} \vert\vert_{ W^{1,p}(\Omega)}  \geq \int_{\Omega} \vert u_{k} \vert^{p}dx > k  \int_{\Omega} \sum_{i=1}^{n} \vert \frac{\partial u_{k}}{\partial x_{i}} \vert^{p}dx }
\end{equation}

Equivalently, this can be written as:\\
\begin{equation}
\displaystyle{ \int_{\Omega} \sum_{i=1}^{n} \vert \frac{\partial u_{k}}{\partial x_{i}} \vert^{p}dx \leq \frac{1}{k} }
\end{equation}

\n Using the result of compactness of the embedding $\displaystyle{W^{1,p}(\Omega) \subset L^{p}(\Omega)}$, we are thus able to find a subsequence $\{u_{k_{l}}\}$ convergent in $L^{p}(\Omega)$ to $u$. Thanks to above equation we have also, $u_{k_{l}} \longrightarrow u$ in $W^{1,p}(\Omega)$. But this same estimate shows also, that $\displaystyle{\frac{\partial u}{\partial x_{i}} = 0 \,\, a.e. \,\text{in}\,\, \Omega}$. It can be shown, that $u = const$: in $\Omega $. Because of the zero mean (the property preserved from the sequence) $u$ must be equal zero (a.e. in $\Omega$). This contradicts the property$\displaystyle{\vert\vert u_{k} \vert\vert_{ W^{1,p}(\Omega)} = 1}$.\\

 \n {\color{red}$\bullet$} We finish this section by given some fundamental properties of Sobolev spaces and we send to {\color{blue}[[1]] } for their proofs.\\

\begin{proposition} ({\color{red} Compact embeddings of Sobolev spaces})\\

\n Let $\Omega$ be a bounded domain in $\mathbb{R}^{n}$ having the cone property. \\

\n  Then the following embeddings are compact:\\

\n (i)  $\displaystyle{W^{j + m,p}(\Omega) \subset W^{j,q}(\Omega)}$ if \quad $0 < n - mp$ and $\displaystyle{j + m - \frac{n}{p} \geq j - \frac{n}{q}}$,\\

\n (ii) $\displaystyle{W^{j + m,p}(\Omega) \subset C^{j}(\overline{\Omega})}$ \quad $mp > n$.\\
\end{proposition}

 \begin{lemma}({\color{red}Gagliardo-Nirenberg inequalities })\\
\begin{equation}
 \displaystyle{\vert\vert u \vert\vert_{L_{\infty}} \leq \vert\vert u \vert\vert^{\frac{1}{2}}\vert\vert Du \vert\vert^{\frac{1}{2}} \,\, \text{for} \,\, u \in H^{s}(\mathbb{R}) ; s \geq 1 \,\, \text{where} \,\, D = \frac{d}{dx}}
 \end{equation}
 \begin{equation}
 \displaystyle{\vert\vert D^{j}u \vert\vert \leq M\vert\vert u \vert\vert^{1-\frac{1}{m}}\vert\vert D^{m}u \vert\vert^{\frac{1}{m}} \,\, \text{for} \,\, u \in H^{m}(\mathbb{R}) ; 0 \leq j \leq m}
 \end{equation}
 \end{lemma}
 \n {\color{red}$\bullet$} (Gagliardo-Nirenberg inequality {\color{blue}[[10]]}): \\

\n Let $\Omega$  be a bounded domain with $\partial\Omega$  in $C^{m}$, and let $u$ be any function in $\displaystyle{W^{m,r}(\Omega) \cap L^{q}(\Omega), 1 \leq q , r \leq \infty}$.\\
\n  For any integer $j, 0 \leq j < m$, and for any number $a$ in the interval $\displaystyle{\frac{j}{m} \leq a \leq 1}$, set\\

$$\displaystyle{\frac{1}{p} = \frac{j}{n} + a(\frac{1}{r} - \frac{m}{n} + (1-a)\frac{1}{q}}$$

\n (i) If $\displaystyle{ m - j - \frac{n}{r} }$ is {\color{red}not} a nonnegative integer, then \\

\begin{equation}
 \displaystyle{\vert\vert D^{j}u \vert\vert_{L^{p}(\Omega)} \leq  C\vert\vert u \vert\vert_{W^{m,r}(\Omega)}^{a} \vert\vert u \vert\vert_{L^{q}(\Omega)}^{1-a}}
\end{equation}

\n (ii) If $\displaystyle{ m - j - \frac{n}{r} }$ is  a nonnegative integer, then the above inequality holds for $\displaystyle{a = \frac{j}{m}}$. The constant $C$ depends only on $\Omega, r, q, j$, and $a$.\\

\n {\color{red}$\bullet$} In the sequel, we will use the following inequalities as the specific cases of the Gagliardo-Nirenberg inequality:\\

\begin{equation}
 \displaystyle{\vert\vert D^{j}u \vert\vert_{L_{\infty}(\Omega)} \leq  C\vert\vert u \vert\vert_{H^{m}(\Omega)}^{a} \vert\vert u \vert\vert_{L^{2}(\Omega)}^{1-a}, \quad ma = j + \frac{n}{2}}
\end{equation}

\begin{equation}
 \displaystyle{\vert\vert D^{j}u \vert\vert_{L^{2}(\Omega)} \leq  C\vert\vert u \vert\vert_{H^{m}(\Omega)}^{a} \vert\vert u \vert\vert_{L^{2}(\Omega)}^{1-a}, \quad ma = j}
\end{equation}

\begin{equation}
 \displaystyle{\vert\vert D^{j}u \vert\vert_{L^{4}(\Omega)} \leq  C\vert\vert u \vert\vert_{H^{m}(\Omega)}^{a} \vert\vert u \vert\vert_{L^{2}(\Omega)}^{1-a}, \quad ma = j + \frac{n}{4}}
\end{equation}

\begin{lemma}({\color{red}Kato's inequality})\\

\n   $\displaystyle{ \text{if} \,\, k > 2 \,\, \text{and}\,\, u \in H^{k}(\mathbb{R})}$ \\
 \begin{equation}
 \displaystyle{\vert < u, uDu >_{k} \vert \leq C \vert\vert u \vert\vert_{2} \vert\vert u \vert\vert_{k}^{2}} \
 \end{equation}
  \n where $\displaystyle{ < . >_{k}}$  denotes the scalar product in  $H^{k}(\mathbb{R})$.\\
  \end{lemma}
  
  \begin{lemma}({\color{red}Friedrich's Inequality})\\
  
 \n  For any $ \epsilon > 0$ there exists an integer $M > 0$ and a real-valued functions $\displaystyle{w_{1}, w_{2}, ....., w_{M} \, \, \text{in}\,\, L^{2}(\Omega)}$, the domain $\Omega$ bounded, that $\displaystyle{\vert\vert w_{j} \vert\vert_{L^{2}(\Omega)} = 1}$, and for any real-valued function $\displaystyle{u \in H_{0}^{1}(\Omega)}$:\\
 \begin{equation}
 \displaystyle{\vert\vert u \vert\vert _{L^{2}(\Omega)}^{2} \leq \epsilon  \vert\vert u \vert\vert _{H^{1}(\Omega)}^{2} + \sum_{j=1}^{M} < u , w_{j} >}
 \end{equation}
 \end{lemma}
 \begin{theorem} ({\color{red}Rellich's Theorem})\\
 
 \n Let $\{u_{m}\}$ be a sequence of functions in $\displaystyle{H_{0}^{1}(\Omega)}$ ($\Omega$ bounded) such that $\displaystyle{\vert\vert u \vert\vert_{H^{1}(\Omega)} \leq const: < \infty}$.\\
There exists a subsequence $\{u_{m'}\}$  convergent in $L^{2}(\Omega)$.
 \end{theorem}
 
 \section{ {\bf{\color{red}  Bonna-Smith results on KDV equation (1975)}}}
 
\n Considered here are waves on the surface of an inviscid fluid in a flat channel. When one is interested in the propagation of one-directional irrotational small amplitude long waves, it is classical to model  the waves by the well-known KdV (Korteweg-de Vries) equation (see {\color{blue}[[38]]}).\\

\n{\color{red}$\bullet$}:  The expression of the Korteweg de Vries equation (KdV) equation is given by:\\ 
 $$\displaystyle{u_{t} + u_{x} + \epsilon uu_{x} + \mu u_{xxx} = 0}$$
 \n where $\displaystyle{u: \mathbb{R}^{2} \longrightarrow \mathbb{R}, (t, x) \longrightarrow u(t, x)}$ is a real function  and $u_{t}$ is its derivative on $t$, $u_{x}$ is its derivative on $x$ and $u_{xxx}$ is its theird derivative on $x$.\\

\begin{lemma}
\n The expression of the KdV equation can be simplified as:\\
$$\displaystyle{u_{t} +  \epsilon uu_{x} + \mu u_{xxx} = 0}$$
\n by means of the substitution: $s = x - t$.\\
\end{lemma}
\n {\color{red}{\bf Proof}}\\

\n {\color{red}$\rhd_{1}$} $\displaystyle{u(t, x) =  v(t, s(t,x)) \Longrightarrow \frac{\partial}{\partial t}u(t, x) = \frac{\partial}{\partial t} [v(t, s(t,x)) ]}$\\

\n $\displaystyle{ = \frac{\partial}{\partial t}v(t, s(t,x)} + \frac{\partial}{\partial t}s(t, x)\frac{\partial}{\partial s(t,x)}v(t,s(t,x))$

\n $\displaystyle{= \frac{\partial}{\partial t}v(t, s) - \frac{\partial}{\partial s}v(t,s)}$.\\

\n \n {\color{red}$\rhd_{2}$} $\displaystyle{u(t, x) =  v(t, s(t,x)) \Longrightarrow \frac{\partial}{\partial x}u(t, x) = \frac{\partial}{\partial x} [v(t, s(t,x)) ] }$\\

\n $\displaystyle{= \frac{\partial t}{\partial x}v(t, s(t,x)) + \frac{\partial}{\partial x}s(t, x)\frac{\partial}{\partial s(t,x)}v(t,s(t,x))}$\\
 
\n  $\displaystyle{ = 0 + \frac{\partial}{\partial s} v(t,s) }$\\

\n i.e.\\

\n  $\displaystyle{\frac{\partial}{\partial t}u(t, x) = \frac{\partial}{\partial t}v(t, s) - \frac{\partial}{\partial s}v(t,s)}$.\\

\n  $\displaystyle{\frac{\partial}{\partial x}u(t, x) = \frac{\partial}{\partial s} v(t,s) }$\\

\n Hence\\

\n  $\displaystyle{ \frac{\partial}{\partial t}v(t, s) - \frac{\partial}{\partial s}v(t,s) + \frac{\partial}{\partial s} v(t,s)  + \epsilon v\frac{\partial}{\partial s}v(t,s) + \mu \frac{\partial^{3}}{\partial s^{3}}v(t,s)}$\\

\n Or\\

\n  $\displaystyle{ \frac{\partial}{\partial t}v(t, s) + \epsilon v\frac{\partial}{\partial s}v(t,s) + \mu \frac{\partial^{3}}{\partial s^{3}}v(t,s)}$\\

\n {\color{red}$\bullet$} Consider the initial-value problem for the forced Korteweg de Vries equation :\\
\begin{equation}
\left \{ \begin{array} {c} \displaystyle{ \frac{\partial}{\partial t}u(t, x) + u\frac{\partial}{\partial x}u(t,s) + \frac{\partial^{3}}{\partial x^{3}}u(t,x) = f(t, x)}\\
\quad\\
u(0, x ) = u_{0}(x) \quad\quad \quad\quad \quad\quad \quad\quad \quad\quad \quad\quad\\
\end{array} \right.
\end{equation}
. \hfill (1.1)\\

\n for $t, x \in \mathbb{R}$. Here the dependent variable $u = u(t, x)$ is a real-valued function of the independent $t$ and $x$ that in most  situations where the equation appears as model, correspond to distance measured in the direction of  the waves propagation and elapsed time. The forcing term $f$ may be thought of as providing a rough accounting of terms that are neglected in arriving at the tidy KDV equation bellow :\\

$$\left \{ \begin{array} {c} \displaystyle{u_{t} + u_{x} +  uu_{x} +  u_{xxx} = 0}\\
\quad\\
u(0, x) = u_{0}(x) \quad\quad  \quad\quad  \quad\\
\end{array} \right.$$

\n The expression of this equation can be simplified by above lemma to following expression:\\
\begin{equation}
\left \{ \begin{array} {c} \displaystyle{u_{t} +  uu_{x} +  u_{xxx} = 0}\\
\quad\\
u(0, x) = u_{0}(x) \quad\quad  \quad \\
\end{array} \right.
\end{equation}

. \hfill { } (1.2)\\

\n by means of the substitution: $s = x - t$.\\

\n The well-posedness of above equation in classical, $L_{2}$-base Sobolev spaces $H^{s}(\mathbb{R})$ for $\displaystyle{s > \frac{3}{2}}$ was  well established  in the mid-1970's (see {\color{blue}[[5]]}, {\color{blue}[[14]]} {\color{blue}[[12]]}  and  the references therein). \\

\n In the early 1980's  Kato {\color{blue}[[13]]} discovered a {\color{red}subtle and  rather general smoothing effect} for above equation.\\

\begin{definition} ( {\color{red}Sharp Kato smoothing property and Kato smoothing property})\\

\n The solutions of the Cauchy problem $u_{t}(t, x) + A u(t, x) = 0$ where $A$  is an operator linear or nonlinear on $[0, T]\times\mathbb{T}$ with $\mathbb{T}$ is a periodic domain, possess :\\

\n (i) the sharp Kato smoothing property if \\
$$\displaystyle{u_{0} \in H^{s}(\mathbb{T}) \Longrightarrow  \partial_{x}^{s+1}u \in L_{loc}^{\infty}(\mathbb{T} ; L^{2}(0 , T)) , T > 0}$$
\n (ii) the Kato smoothing property if\\
 $$\displaystyle{u_{0} \in H^{s}(\mathbb{T}) \Longrightarrow  u \in L^{2}(0, T, H^{s+1}(\mathbb{T}))}$$
\end{definition}

\begin{remark}
It is well known that\\

\n (i) the solutions of the Cauchy problem of the KdV equation on a periodic domain $\mathbb{T}$:
$$\displaystyle{u_{t} + uu_{x} + u_{xxx} = 0, \quad u(0, x) = u_{0}(x), \quad x \in \mathbb{T} , t \in \mathbb{R}}$$ 
possess neither the sharp Kato smoothing property nor the Kato smoothing property.\\
\n (ii) but for the solutions of the Cauchy problem of the KdV-Burgers (KdVB) equation on a periodic domain $\mathbb{T}$:\\
$$\displaystyle{u_{t} + uu_{x} + u_{xxx} - u_{xx} = 0, \quad u(0, x) = u_{0}(x), \quad x \in \mathbb{T} , t \in \mathbb{R}}$$ 
possess the sharp Kato smoothing property and the Kato smoothing property\\
\n We present (ii) in following theorem:\\
\end{remark}

\n Let $x \in [0, L] , L > 0$ then $\displaystyle{\frac{x}{L} \in [0, 1]}$ and let $\displaystyle{e_{n}(x) = e^{2in\pi x}, n = 0, \pm 1, \pm 2, .....}$ with $x \in [0, 1]$ Then $\displaystyle{\{e_{n}(x)\}_{-\infty}^{\infty}}$ forms an orthonormal basis in the space $L^{2}(0; 1)$.We may define the Sobolev space $\displaystyle{H_{p}^{s} : = H_{p}^{s}(0, 1)}$ of order $s$ ($s \in \mathbb{R}$ ) as the space of all real periodic functions of period $1$.\\
$$\displaystyle{v(x) = \sum_{n = -\infty}^{\infty}v_{n} e_{n}(x)}$$
\n such that\\
\begin{equation}
\displaystyle{\vert\vert v \vert\vert_{s} : = \sum_{n = -\infty}^{\infty} ( 1 + \vert n \vert ^{2s})^{\frac{s}{2}} \vert v_{n}\vert^{2} < \infty}
\end{equation}
$\displaystyle{\vert\vert v \vert\vert_{s}}$ is a Hilbert norm for $H_{p}^{s}$.\\

\n For any $\displaystyle{s \in \mathbb{R},  v \in H_{p}^{s}}$ with \\
$$\displaystyle{v(x) = \sum_{k = -\infty}^{\infty}v_{k} e_{k}(x)}$$

\n we define the operator $D_{x}^{s}$ by\\
$$\displaystyle{D_{x}^{s}v(x) = \sum_{k = -\infty}^{\infty}(1 + k^{2})^{s/2}v_{k} e_{k}(x)}$$
\n such that:\\
$$\displaystyle{\vert\vert v \vert\vert_{s} = \vert\vert D_{x}^{s}v \vert\vert_{L^{2}(0, 1)}}$$
\begin{theorem}

\n Consider the Cauchy problem of the KdV-Burgers equation posed on $(0; L), L > 0$ with periodic boundary conditions,\\
\begin{equation}
\left \{ \begin{array} {c} \displaystyle{u_{t} + uu_{x} + u_{xxx} - u_{xx} = 0, \quad x \in (0, L)}  \quad\quad  \quad\quad \quad\quad \quad\\
\quad\\
u(0, x) = u_{0}(x) \quad\quad \quad\quad  \quad\quad  \quad\quad  \quad\quad  \quad\quad  \quad\quad \quad\quad \quad\\
\quad\\
u(t; 0) = u(t; L); u_{x}(t; L) = u_{x}(t; L); u_{xx}(t; 0) = u_{xx}(t; L)\\
\end{array} \right.\\
\end{equation}
\n which is well-known to be globally well-posed in the space $H_{p}^{s}$  for $s \geq 0$:
Let $s\geq 0$ and $T > 0$ be given. Then\\

\n (i) For any $u_{0} \in  H_{p}^{s}$ the Cauchy problem of the KdV-Burgers equation admits a unique mild solution\\
$$\displaystyle{u \in C([0, T],  H_{p}^{s})}$$
\n which possesses the Kato smoothing property:\\
$$\displaystyle{u_{0} \in H_{p}^{s} \Longrightarrow u \in L^{2}(0, T;  H_{p}^{s + 1} )}$$
\n It also possesses the sharp Kato smoothing property.\\

\n (ii) For any $u_{0} \in  H_{p}^{s}$ the corresponding solution $u$ of Cauchy problem of the KdV-Burgers equation belongs to the space $\displaystyle{C([0, T];  H_{p}^{s}) \cap L^{2}([0, T];  H_{p}^{s + 1})}$ and satisfies \\
$$\displaystyle{sup_{0 \leq x \leq L} \vert\vert \partial_{x}^{s+1}u(., x) \vert\vert_{L^{2}(0, T)} \leq \omega_{s, T} \vert\vert u_{0} \vert\vert_{s}}$$
\n where  $\displaystyle{\omega_{s;T} : \mathbb{R}_{+} \longrightarrow \mathbb{R}_{+} }$ is a nondecreasing continuous function.\\

\end{theorem}

\n By contrast, the theory pertaining  to the (1.1) for the forced KDV equation has remained less developed. The following result was given by {\color{blue}[[5]]} in the early 1970's:\\

\begin{theorem}

\n For given $T > 0$ and $s \geq 3$, if \\

\n (i) $\displaystyle{u_{0} \in H^{s}(\mathbb{R})}$,\\

\n (ii) $\displaystyle{f \in C(- T, T ; H^{s}(\mathbb{R}))}$,\\

\n (iii) $\displaystyle{f_{t} \in C(- T, T ; L^{2}(\mathbb{R}))}$,\\

\n then\\

\n The initial-value problem (1.1) has a unique solution  $\displaystyle{u \in C(- T, T ; H^{s}(\mathbb{R})) \cap C(- T, T ; L^{2}(\mathbb{R}))}$.\\

\n In addition, the solution $u$ depends continuously in $C(- T, T ; H^{s}(\mathbb{R})$ on $u_{0}$ in $H^{s}(\mathbb{R})$ and $f$ in $\displaystyle{u \in C(- T, T ; H^{s}(\mathbb{R})) \cap C^{1}(- T, T ; L^{2}(\mathbb{R}))}$.\\
\end{theorem}
\begin{remark}
\n This result was strengthened by {\color{blue}[[39]]} where he showed that the conclusion of this theorem holds without assumption (iii).\\
\end{remark}

\n {\color{red}$\rhd$} It  frequently arises in numerical simulation where  unbounded domains are hard  to model and the relative  simplicity of imposing periodic boundary conditions is very attractive.\\

\n Let $H^{s}(\mathbb{S})$ denote the real Sobolev space of order $s$ ($s \geq 0$) on the unit length circle in the plane. $H^{s}(\mathbb{S})$ may be  characterized as the space of real $1$-periodic functions
 $v$ whose Fourier series\\
 $$ \displaystyle{v(x) \sim \sum_{-\infty}^{\infty} v_{k}e^{2i\pi kx}}$$
  \n is such that\\
  $$ \displaystyle{\vert\vert v\vert\vert_{s} = \{\sum_{-\infty}^{\infty} (1 + \vert k\vert^{2s}\vert v_{k}\vert^{2}\}^{\frac{1}{2}}}$$
  \n  $\displaystyle{\vert\vert v\vert\vert_{s}}$ defines a Hilbert space norm on the linear space $H^{s}(\mathbb{S})$ .\\
 
 \n Let $D^{s}$ represent the fractional derivative of order $s$, so if $v$ has the Fourier series as above then\\
  $$ \displaystyle{D^{s}v(x) \sim \sum_{-\infty}^{\infty} v_{k}\vert k \vert^{s}e^{2i\pi kx}}$$
 \n By using energy estimates we have the following theorem\\

\begin{theorem}

\n For $s \geq 2$ and $T > 0$, if $ u_{0} \in H^{s}(\mathbb{S})$ and $f \in L^{1}(-T, T; H^{s}(\mathbb{S}))$, then the initial-value problem :\\
\begin{equation}
\left \{ \begin{array} {c} \displaystyle{u_{t} + uu_{x} + u_{xxx} = f \quad t \in \mathbb{R}, x \in \mathbb{S}}\\
\quad\\
u(0,x) = u_{0}(x) \quad\quad \quad\quad \quad\quad \quad\quad \quad\\
\end{array} \right.
\end{equation}
\end{theorem}
\n has a unique solution  $u \in  C(-T, T; H^{s}(\mathbb{S}))$. More over, the solution $u$ depends continuously on $u_{0}$ in $H^{s}(\mathbb{S})$ and $  \in L^{1}(-T, T; H^{s}(\mathbb{S}))$\\

\section{{\bf{\color{red}  Aimar's results on K.S equation (1982)}}}
  
  \n{\color{red}$\rhd_{1}$} The equation\\
  
  $\displaystyle{u_{t} + \alpha u_{xx} + \beta u_{xxxx} + \gamma uu_{x} = 0}$ \hfill { } {\color{blue} (KS)}\\
  
  \n arises in interesting physical situations, for example as a model for long waves on a viscous fluid flowing down an inclined plane  {\color{blue}[[37]]} and to derive drift waves in a plasma  {\color{blue}[[6]]}. Equation {\color{blue} (KS)} was derived also independently by Kuramoto  {\color{blue}[[22]]} and Kuramoto-Tsuzuki {\color{blue}[[23]]}, {\color{blue}[[24]]} as a model for phase turbulence in reaction-diffusion systems and by Sivashinsky  {\color{blue}[[33]]} as a model for plane flame propagation, describing the combined influence of diffusion and thermal conduction of the gas on the stability of a plane flame front.\\
  
\n{\color{red}$\rhd_{2}$} The Michelson system arises as an equation for the derivative of a traveling wave or a steady state solution of the one-dimensional Kuramoto-Sivashinsky equation (see {\color{blue}[[2]]}, {\color{blue}[[3]]}, {\color{blue}[[21]]}, {\color{blue}[[28]]}, {\color{blue}[[32]]}, {\color{blue}[[33]]} and the references therein).\\

\begin{equation}
 \displaystyle{ u_{t} = - u_{xxxx} - u_{xx} - uu_{x}}
\end{equation}

\n This one-dimensional Kuramoto-Sivashinsky (K-S) equation has been rigorous studied of view mathematic (existence and uniqueness  of the solution on some Sobolev spaces) by Aimar in its PhD (1982) under the following form:\\

\begin{theorem} ({\bf{\color{red}Aimar 1982}})\\

\n Let
\begin{equation}
\displaystyle{ u_{t} = - u_{xxxx} - u_{xx} - \frac{1}{2}[u'(x)]^{2} ; (t,x) \in [0, T]\times[0, L]}
\end{equation}

\n be the one-dimensional Kuramoto-Sivashinsky  equation with\\

\n $u_{x}(t,0) = u_{x}(t,L) = 0$ and $u_{xxx}(t, 0) = u_{xxx}(t, L) = 0$  as limit conditions and $T, L$ are parameters.\\

\n Let $\mathbb{H}^{s} := \mathbb{H}^{s}(0, L)$ be the classical Sobolev space where $\mathbb{H}^{0}:= \mathbb{L}^{2}$.\\

\n Let $\displaystyle{\mathbb{V} = \{ v \in \mathbb{H}^{1}; v_{xxx}\in \mathbb{L}^{2}\}}$\\

\n If we choose $u(0, x)$ be in $\mathbb{H}^{1}$\\

\n Then there exists an uniqueness solution of the problem {\color{red}(4.2)} such that\\

\n i) $\displaystyle{u \in \mathbb{L}^{2}(0, T; \mathbb{V})\cap \mathbb{L}^{\infty}(0, T; \mathbb{H}^{2})}$\\

\n ii) $\displaystyle{u_{t} \in \mathbb{L}^{2}([0, T]\times [0, L])}$\\
\end{theorem}
\begin{remark} 
\quad\\

\n {\color{red}(i)} As $\displaystyle{u \in \mathbb{L}^{\infty}(0, T; \mathbb{L}^{2})}$ then $\displaystyle{u : [0, T] \longrightarrow \mathbb{L}^{2}}$ is continuous.\\

\n{\color{red}(ii)} If $\displaystyle{u \in \mathbb{L}^{\infty}(0, T; \mathbb{H}^{2})}$ then $\displaystyle{u : [0, T] \longrightarrow \mathbb{H}^{1}}$ is continuous and the choose of $u(0, x)$ in $\mathbb{H}^{1}$ is justified.\\

\n{\color{red}(iii)} If we replace the limit conditions of above theorem by:\\

\n $\left\{\begin{array}[c]{l}\displaystyle{u(t,0) = u(t, L)}\\
\quad\\
\displaystyle{u_{x}(t,0) = u_{x}(t, L)}\\
\quad\\
\displaystyle{u_{xx}(t,0) = u_{xx}(t, L)}\\
\quad\\
\displaystyle{u_{xxx}(t,0) = u_{xxx}(t, L)}\\
\end{array}\right.$ \\

\n Then we have existence and uniqueness  of the solution of kuramoto-Sivashinsky with these limit conditions.\\

\n{\color{red}(iv)} In {\color{red}(4.2)} the second-order term acts as an energy source and has a destabilizing effect, the nonlinear term transfers energy from low to high wave numbers, while the  fourth-order term removes the energy on small scales. The Kuramoto-Sivashinsky equation is known to exhibit spatio-temporal chaos.\\

\n By adding a periodic second member $f$ in {\color{red}(4.2)} and taking $L=1$ and the following conditions \\

\n $\left\{\begin{array}[c]{l}\displaystyle{u(t,0) = u(t, 1)}\\
\quad\\
\displaystyle{u_{x}(t,0) = u_{x}(t, 1)}\\
\quad\\
\displaystyle{u_{xx}(t,0) = u_{xx}(t, 1)}\\
\quad\\
\displaystyle{u_{xxx}(t,0) = u_{xxx}(t, 1)}\\
\end{array}\right.$ \\
\end{remark}

\n {\color{red}($\rhd$} Now, we consider numerical solutions of Kuramoto-Shivasinsky equation, we will give a semi-implicit numerical method with a result of stability. Namely, the numerical approximation is taken on regular grid by using an implicit finite differences for linear operator and semi-implicit finite differences for non linear operator.\\
 
\n Let $h = \frac{1}{N}$, $\displaystyle{\Omega = ]0, 1[}$ and $\displaystyle{\Omega_{h} = \{x_{p}; x_{p} = (p-1)h, p = 1, ..., N+1\}}$ \\

\n Let $\sigma_{h}(x_{p})$ be the interval of center $x_{p}$ and of diameter $h$.\\

\n Let $w_{h}^{p}$ be the characteristic function of $\sigma_{h}(x_{p})$ defined by:\\

$\left\{\begin{array}[c]{l}\displaystyle{w_{h}^{p}(x) = 1 \quad if \quad x \in [x_{p} - \frac{1}{2}, x_{p} + \frac{1}{2}[}\\
\quad\\
\displaystyle{w_{h}^{p}(x) = 0 \quad if \quad x /\!\!\!\!\!\in [x_{p} - \frac{1}{2}, x_{p} + \frac{1}{2}[}\\
\end{array}\right.$ \\

\n We put $\displaystyle{\mathbb{V}_{h} = span\{w_{h}^{1}, w_{h}^{2}, ...,w_{h}^{N+1}\}}$ and we define $u_{h}(t)$ in $\mathbb{V}_{h}$ by $\displaystyle{u_{h}(t) = \sum_{p=1}^{N+1}u_{h}^{p}(t)w_{h}^{p}}$\\

\n The usual scalar product on $\mathbb{V}_{h}$ is given by\\

\n $\displaystyle{ <u_{h}(t), u_{h}(t)> = h\sum_{p=1}^{N}u_{h}^{p}(t)v_{h}^{p}(t)}$\\

\n with associated norm:\\

\n $\displaystyle{||u_{h}(t)||_{h}^{2} = h\sum_{p=1}^{N}|u_{h}^{p}(t)|^{2}}$\\

\n $||.||_{h}$ is an approximation of the usual norm in $\mathbb{L}^{2}(\Omega)$.\\

\n For $u_{h}(t)$ in $\mathbb{V}_{h}$, we define following classical operators:\\

\n $\left\{\begin{array}[c]{l}\displaystyle{\delta_{h}^{+}u_{h}^{p}(t) = \frac{u_{h}^{p+1}(t) - u_{h}^{p}(t)}{h}}\\
\quad\\
\displaystyle{\delta_{h}^{-}u_{h}^{p}(t) = \frac{u_{h}^{p}(t) - u_{h}^{p-1}(t)}{h}}\\
\quad\\
\displaystyle{\Delta_{h}u_{h}^{p}(t) = \frac{u_{h}^{p+1}(t) - 2u_{h}^{p}(t) + u_{h}^{p-1}(t)}{h^{2}}}\\
\quad\\
\displaystyle{\Delta_{h}^{2}u_{h}^{p}(t) = \frac{u_{h}^{p+2}(t) - 4u_{h}^{p+1}(t) + 6u_{h}^{p}(t) - 4u_{h}^{p-1}(t) + u_{h}^{p-2}(t)}{h^{4}}}\\
\end{array}\right.$ \\

\n We turn now to consider the nonlinear term of Kuramoto-Shivasinsky equation, we choose the following semi-implicit finite differences :\\

\n {\color{red}$\bullet$} $\displaystyle{[D_{h}u_{h}^{p}(t)]^{2} = \frac{1}{3}[[\delta_{h}^{+}u_{h}^{p}(t)]^{2} + [\delta_{h}^{-}u_{h}^{p}(t)]^{2} + \delta_{h}^{+}u_{h}^{p}(t)\delta_{h}^{+}u_{h}^{p}(t)]}$\\

\n We observe that $\displaystyle{||\delta_{h}^{+}u_{h}(t)||_{h}^{2} = ||\delta_{h}^{-}u_{h}(t)||_{h}^{2}}$ that we denote by $\displaystyle{||\delta_{h}^{1}u_{h}(t)||_{h}^{2}}$\\

\n and\\

\n $\displaystyle{||\delta_{h}^{+}\Delta_{h}u_{h}(t)||_{h}^{2} = ||\Delta_{h}\delta_{h}^{-}u_{h}(t)||_{h}^{2}}$ that we denote by $\displaystyle{||\delta_{h}^{1}\Delta_{h}u_{h}(t)||_{h}^{2}}$\\

\n Now we define two norms on $\mathbb{V}_{h}$ as follows:\\

\n $\displaystyle{||u_{h}(t)||_{h,1}^{2} = ||u_{h}(t)||_{h}^{2} + ||\delta_{h}^{1}u_{h}(t)||_{h}^{2}}$ this norm is an approximation of usual norm of $\mathbb{H}^{1}(\Omega)$\\

\n and\\

\n $\displaystyle{||u_{h}(t)||_{h,1}^{2} = ||u_{h}(t)||_{h}^{2} + ||\Delta_{h}u_{h}(t)||_{h}^{2}}$ this norm is an approximation of usual norm of $\mathbb{H}^{2}(\Omega)$\\

\n Then the approximed problem can be written as follows:\\

\n $\left\{\begin{array}[c]{l}\displaystyle{\frac{d}{dt}u_{h}^{p} + \Delta_{h}^{2}u_{h}^{p}(t) + \Delta_{h}u_{h}^{p}(t) + \frac{1}{2}[\delta_{h}u_{h}^{p}(t)]^{2} = f_{h}^{p}(t)}\\
\quad\\
with \quad boundary \quad conditions\\
\quad\\
\displaystyle{u_{h}^{N+1}(t) = u_{h}^{1}(t)}\\
\quad\\
\displaystyle{u_{h}^{N+2}(t) = u_{h}^{2}(t)}\\
\quad\\
\displaystyle{u_{h}^{0}(t) = u_{h}^{N}(t)}\\
\quad\\
\displaystyle{u_{h}^{-1}(t) = u_{h}^{N-1}(t)}\\
\quad\\
with \quad initial \quad conditions\\
\quad\\
\displaystyle{u_{h}^{p}(0) = u_{0,h}^{p}; p=1, 2, ...., N}\\
\quad\\
\end{array}\right.$ $\hfill { } {\color{red} (\mathbb{P}_{d})}$\\

\n Let $f_{h}^{*}, u_{0,h}^{*}$ are a perturbation of second member and the initial condition of the equation respectively.\\

\n If we put $\displaystyle{w_{h} = u_{h} - u_{h}^{*}}$ and $\displaystyle{g_{h} = f_{h} - f_{h}^{*}}$, we get the following result of the stability:\\

\begin{theorem} ({\color{blue}Aimar 1982})\\

\n If $\displaystyle{||u_{0,h}(t)||_{h,2}}$, $\displaystyle{||u_{0,h}^{*}(t)||_{h,2}}$, $\displaystyle{\int_{0}^{T}||f_{h}(t)||_{h,2}dt}$ and $\displaystyle{\int_{0}^{T}||f_{h}^{*}(t)||_{h,2}dt}$ are uniformly bounded. Then\\

\n (i) $u_{h}$ and $u_{h}^{*}$ are defined on $[0, T]$\\

\n (ii) $\displaystyle{\forall \epsilon > 0 \quad \exists \eta > 0 }$ such  that\\

\n  $\displaystyle{ \forall w_{0,h}}$  and $\displaystyle{\forall g_{h}}$ satisfying $\displaystyle{ max (||w_{0,h}||_{h,1}^{2} , \int_{0}^{T}||g_{h}(\tau)||_{h,1}^{2}d\tau) < \eta}$\\

\n  implies \\

\n  $\displaystyle{||w_{h}(t)||_{h,1}^{2} \leq \epsilon ; \forall t \in [0, T], \forall h \in ]0, h_{0}[, h_{0} > 0}$\\
\end{theorem}

\n  {\bf {\color{red}The essential points of the proof}}\\

\n  Under the following assumptions \\

 \n {\color{blue}(a)} $\displaystyle{||u_{0,h}(t)||_{h,2} \leq C}$ ; \\
 \n $C$ is independent of $h$ (similar for $u_{0,h}^{*}$) \\

\n  {\color{blue}(b)} $\displaystyle{\int_{0}^{T}||f_{h}^{*}(t)||_{h,2}dt \leq C}$; \\
\n $C$ is independent of $h$ (similar for $f^{*}$) \\

 \n we have\\

\n  {\color{red}($\alpha$)} $\displaystyle{\forall t \in (0, T) , T < +\infty; ||D_{h}^{2}u_{h}(t)||^{2} \leq C}$; \\
\n $C$ is independent of $h$\\

 \n From {\color{red}($\alpha$)}, we deduce the following inequalities\\

\n  {\color{red}($\alpha_{1}$)} $\displaystyle{||[D_{h}u_{h}(t)]^{2}-[D_{h}u_{h}^{*}(t)]^{2}||^{2} \leq C ||\delta_{h}^{1}(u_{h}(t)-u_{h}^{*}(t))||_{h}^{2}}$; \\
\n $C$ is independent of $h$.\\

\n {\color{red}($\alpha_{2}$)} $\displaystyle{ _{_{_{_{_{_{p\in\{1,...,N\}}}}}}}\!\!\!\!\!\!\!\!\!\!\!\!\!\!\!\!\!\!||\delta_{h}^{1}u_{h}^{p}(t)|| \leq C_{1}}$; \\
\n $C_{1}$ is independent of $h$.\\

\n {\color{red}($\alpha_{3}$)} $\displaystyle{[\delta_{h}^{1}u_{h}^{p}(t)]^{2} \leq C_{2} + ||D_{h}^{2}u_{h}(t)||^{2}\quad \forall p\in \{1, ...,N\}}$; \\
\n $C_{2}$ is independent of $h$.\\

\n \begin{center}

{\color{red}Transition phenomen to under harmonic} (Aimar 1982)\\
\includegraphics[scale=0.085]{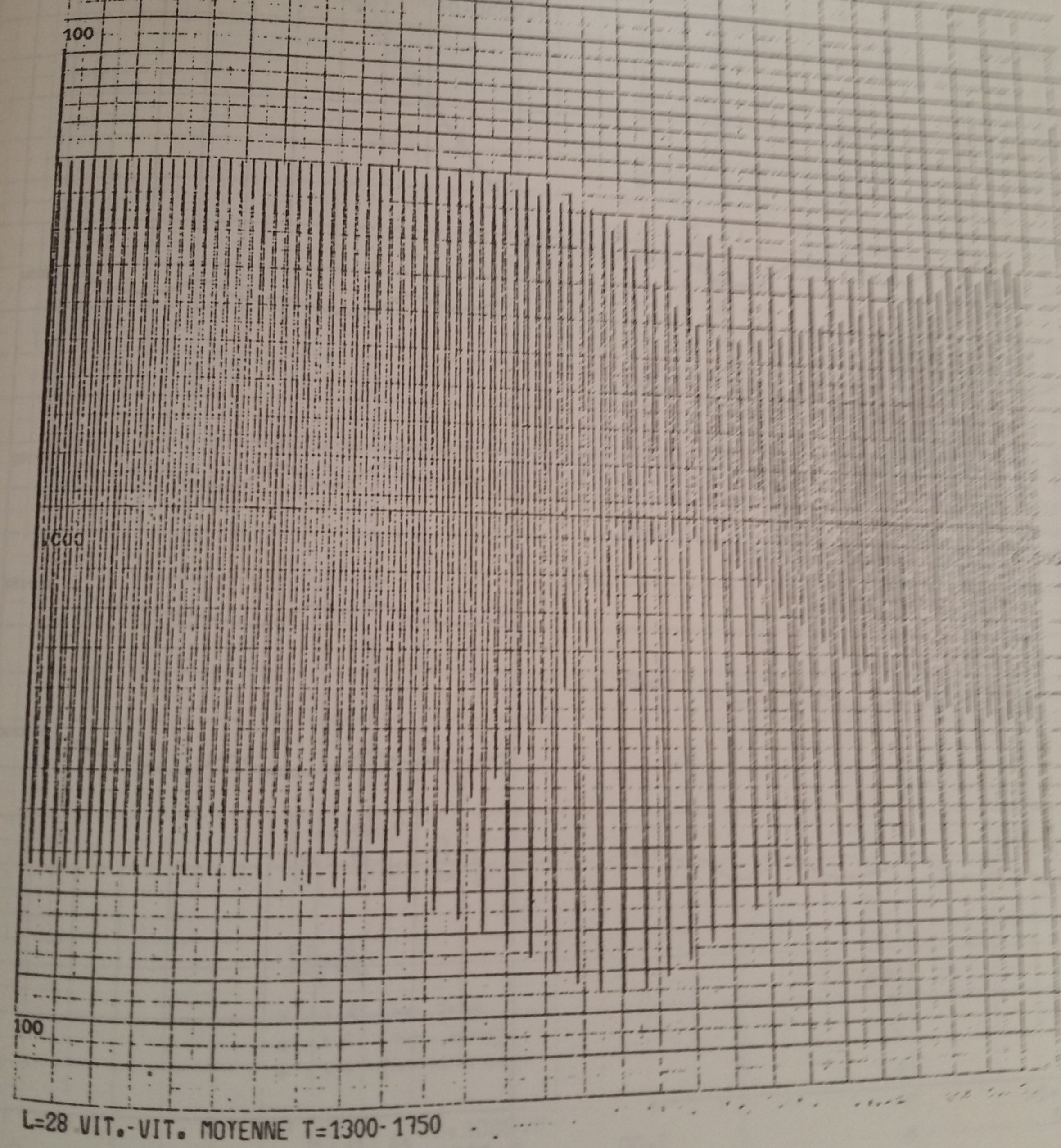}\\
\end{center}

 \begin{center}
 {\color{red} Evolution of flame front} (Aimar 1982)\\

\includegraphics[scale=0.1]{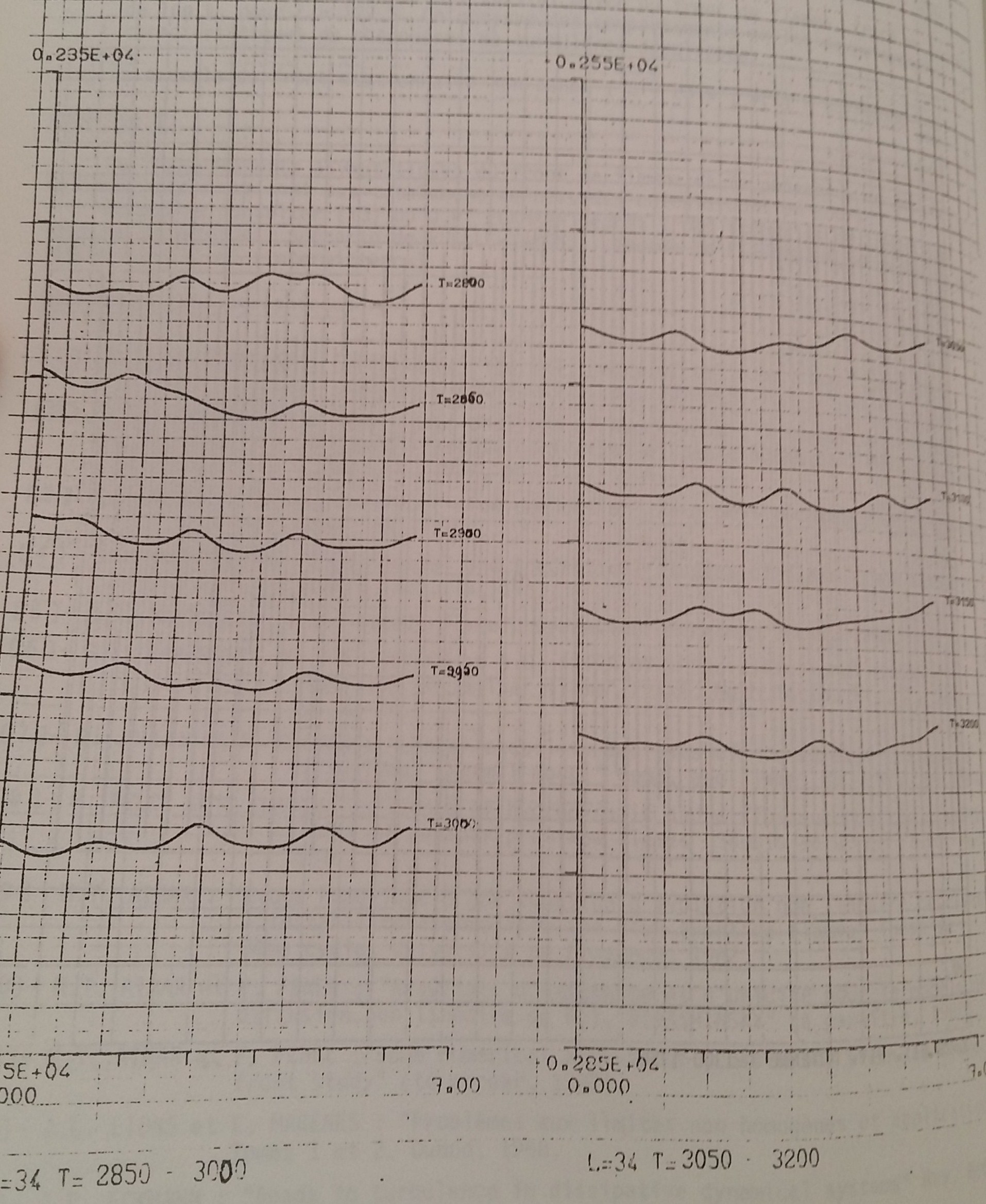}\\
\end{center}
\quad\\

\begin{remark}
\quad\\

\n  We can apply to {\color{red} (4.2)} an another numerical method by using a lower order term $\lambda v_{xx}$, with $\lambda$ chosen in order to counteract the effect of $-v_{xxxx}$.\\

\n  Eq. {\color{red} (4.2)}  is approximed on regular grid, using centered finite differences:\\

\n  $\displaystyle{\frac{v_{j}^{n+1} - v_{j}^{n}}{\delta t} = - v_{j}^{n}\frac{v_{j+1}^{n} - v_{j-1}^{n}}{2\delta x} - \frac{v_{j-1}^{n} - 2v_{j}^{n} + v_{j+1}^{n}}{\delta x^{2}}
  - \frac{v_{j-2}^{n} - 4v_{j-1}^{n} + 6v_{j}^{n} - 4v_{j+1}^{n}+ v_{j+2}^{n}}{\delta x^{4}}}$\\

\n  $\displaystyle{- \lambda\frac{v_{j-1}^{n} - 2v_{j}^{n} + v_{j+1}^{n}}{\delta x^{2}} -\lambda  \frac{v_{j-1}^{n+1} - 2v_{j}^{n+1} + v_{j+1}^{n+1}}{\delta x^{2}}}$\\

\n where $\lambda$ has to be chosen such that the above algorithm is stable.\\

 \n  We see that while this is certainly not the method of choice to solve {\color{red}(4.2)}, it sufficiently accurate to represent the statistics of the solution.\\
 \end{remark} 
 \section{{\color{red}{\bf Presentation of KdV-KS equation }}{\color{blue}{\bf(1996)}}}

\n The above theorem was generalized by Aimar-Intissar (1996) to the Korteweg-de Vries-Kuramoto-Sivashinsky equation (KdV-KS):\\

\n $\displaystyle{ u_{t} = - u_{xxxx} - u_{xxx} - Ru_{xx} - uu_{x}}$\hfill { } {\color{red} (KdV-KS)}\\

\n where $R > 0$ as the ``anti-diffussion'' parameter corresponding to the Reynolds number.\\

\n We give bellow  the theorem of existence and uniqueness of the solution of KdV-KS equation and its proof ($\mathbb{F}$-Method).\\

\n This last equation is an approximate equation for surface waves two dimentional incompressible viscous fluid down an inclined plane under the assumption of small amplitude and long wave. It has a mixed form of Korteweg-de Vries equation (KdV) for surface waves of water.\\

\n $\displaystyle{u_{t} = - u_{xxx} - uu_{x}}$ \hfill { } {\color{red} (KdV)}\\

\n and Kuramoto-Sivashinsky equation {\color{red} (4.1)} for thermodiffusive instability of premixed flame fronts and wave propagation in reaction diffusion systems.\\

\n  Let $L = 2\pi$ and $\mathfrak{H} := \mathbb{L}^{2}(0, 2\pi)$, if we consider the operator $\mathbb{A}u = u_{xxxx} + u_{xxx} + u_{xx}$ with domain $\displaystyle{D(\mathbb{A}) = \{ u \in \mathbb{H}_{loc}^{4}(\mathbb{R}) ; u(x + 2\pi) = u(x) \,\\, \forall \,\, x \in \mathbb{R}\}}$ then we see that it is nonselfadjoint unbounded linear operator with compact resolvent and its real and imaginary parts are linear operators associated to Kuramoto-Shivashinsky equation and Korteweg-de Vries equation respectively.\\

\n The derivative of above equation give {\color{red} (4.1)}. Hence equation {\color{red} (4.1)} has been studied extensively in the effort to understand the complicated dynamics in nonlinear partial differential equations (PDEs).\\

\n Physically, the two forms {\color{red} (4.1)} and {\color{red}(4.2)} of the K-S equation model different things; equation {\color{red} (4.1)} models small disturbances in liquid films falling down an inclined or vertical plane and propagation of concentration waves in  chemical reactions , where equation ({\color{red} (4.2)} models instabilities in laminar flame fronts {\color{blue} [[32]]}, {\color{blue}[[33]]}.\\

\n From equation {\color{red}(4.2)} Aimar-Michelson-Sivashinsky was observed that the mean value of the solution $\displaystyle{v(t) := \frac{1}{2L}\int_{0}^{L}[u(t,x)]dx}$ satisfies the drift equation\\
\begin{equation}
\displaystyle{\frac{d}{dt}v(t) = -\frac{1}{2L}\int_{0}^{L}[u_{x}(t,x)]^{2}dx}
\end{equation}

\n and that solutions took the form\\

\n $u(t,x) = c_{0}^{2} + \phi(t, x)$ where $c_{0}$ is a constant independent of the initial condition,\\

\n and  $\displaystyle{\psi(t):= \frac{1}{2L}\int_{0}^{L}[\phi(t,x)]dx}$ is close to zero.\\

\n Therefore $\displaystyle{\frac{d}{dt}\psi(t) \sim -c_{0}^{2}}$ and hence we
substitute a solution of the form \\ $\displaystyle{u(t, x) = -c^{2}t + \phi(x); c \in \mathbb{R}}$,  we get \\
\begin{equation}
\displaystyle{-c^{2} + \phi_{xxxx} + \phi_{xx} + \frac{1}{2}[\phi_{x}]^{2} = 0}
\end{equation}

\n The energy $E(t)$ of solutions of {\color{red} ( (KdV-KS))} is defined by\\
\begin{equation}
\displaystyle{E(t) := \frac{1}{2L}\int_{0}^{L}[u(t,x)]^{2}dx}
\end{equation}
\n and its higher order energy $\hat{E} (t)$ is defined by\\
\begin{equation}
\displaystyle{\hat{E}(t) := \frac{1}{2L}\int_{0}^{L}[u_{xx}(t,x)]^{2}dx}
\end{equation}
\n The stability of {\color{red} ( (KdV-KS))} significantly depends on the anti-diffusion parameter $R > 0$.\\
\n It was observed in numerical simulations of the KS equation that an orbit will shadow an unstable periodic orbit for a time before diverging from it as some observed dynamics in low-dimensional systems such as the Lorenz equations.\\

\section{{\color{red}{\bf On travelling wave solution of the Korteweg-de Vries-Kuramoto-Sivashinsky equation}}}

\n In this way it is useful to study the simple behaviours of {\color{red} (4.1)}, such as the travelling wave solutions of the form $\tilde{u} = u(x - \omega t)$, and the steady-state solutions.\\

\n We begin by substituting $\tilde{u} = u(\xi)$  where $\xi = x - \omega t$. Then equation {\color{red} (4.1)}
transforms to\\
\begin{equation}
\displaystyle{ - \omega u_{\xi} + u_{\xi\xi\xi\xi} + u_{\xi\xi} + uu_{\xi} = 0}
\end{equation}

\n This can be rewritten as a one-dimensional system by the change of variables $x_{1} = \tilde{u}$, $\displaystyle{x_{2} = \tilde{u}_{\xi}}$, $\displaystyle{x_{2} = \tilde{u}_{\xi\xi}}$ and $\displaystyle{x_{4} = \tilde{u}_{\xi\xi\xi}}$\\

\n Then\\
\begin{equation}
\left\{\begin{array}[c]{l}\displaystyle{\frac{d}{d\xi}x_{1} = x_{2}
}\\
\quad\\
\displaystyle{\frac{d}{d\xi}x_{2} = x_{3} }\\
\quad\\
\displaystyle{\frac{d}{d\xi}x_{3} = x_{4}}\\
\quad\\
\displaystyle{\frac{d}{d\xi}x_{4} = \omega x_{2} - x_{1}x_{2} - x_{3}}\\
\end{array}\right.
\end{equation}

\n Integrating the last equation once gives\\

\n $\displaystyle{\frac{d}{d\xi}x_{3} = \omega x_{1} -  x_{2} - \frac{1}{2}x_{1}^{2}}$ + constant\\

\n We remove the constant by taking $x_{1}, x_{2}, x_{3}, x_{4} = 0$ to be a valid solution for all
$\omega$. Now under the transformation $x = x_{1} - \omega , y = x_{2} , z = x_{3}, t = \xi$ we obtain the system\\
\begin{equation}
\left\{\begin{array}[c]{l}\displaystyle{^{^{\bullet}}\!\!\!x = y}\\
\quad\\
\displaystyle{^{^{\bullet}}\!\!\!y = z}\\
\quad\\
\displaystyle{^{^{\bullet}}\!\!\!z = \frac{\omega^{2}}{2} - \frac{1}{2}x^{2} - y}\\
\end{array}\right.
\end{equation}

\n where the dot denotes always differentiation with respect to time. Thus we have derived the Michelson system with $\omega = c\sqrt{2}$.\\

\n  We will consider equation {\color{red} (4.1)} on periodic boundary conditions, $u(t, x) = u(t, x + L)$ (with $L$ as a parameter). There is a rigorous proof of an inertial manifold {\color{blue}[[31]]}, {\color{blue}[[36]]}  for equation {\color{red} (4.1)} which makes this equation effectively equivalent to a finite dimensional
system for long time behavior.\\

\n Estimates of the dimension of the latter are discussed in {\color{blue}[[36]]}. This, together with the fact that it is a scalar equation, makes it a paradigm for the study of rich spatiotemporal dynamics in one dimensional PDEs, and thus the K-S equation has been well-studied numerically. In particular in {\color{blue}[[2]]}, {\color{blue}[[3]]}.\\
 \n  reported a rich and diverse variety in the nature of the solutions to the K-S equation on periodic boundary conditions, as the control parameter $L$ is varied. \\

  \section{{\color{red}{\bf Existence and Uniquness of KdV-KS equation  on domain $[0, T]\times [0, 2\pi] , T > 0$}}}   
  \n {\color{blue} {\bf $\mathbb{F}$-Method }}\\

\n Integrating KdV-KS equation once, we deduce \\

\n $\displaystyle{ u_{t} + u_{xxxx} + u_{xxx} + Ru_{xx} + \frac{1}{2}u_{x}^{2} = c_{0}}$ $\hfill { } {\color{red}(\mathbb{F}_{1})}$ \\

\n where $c_{0}$ is the integration constant.\\

\n If we introduce $f(t,x)$ in place of $c_{0}$ and $\displaystyle{u(t,x) = e^{\lambda t}v(t,x)}$ where $\lambda \in \mathbb{R}$ we get :\\

\n $\displaystyle{ v_{t} + v_{xxxx} + v_{xxx} + Rv_{xx} + \frac{1}{2}e^{\lambda t}v_{x}^{2} = e^{-\lambda t}f}$ $\hfill { } {\color{red}(\mathbb{F}_{2})}$ \\

\n Let $\displaystyle{\tilde{f} = e^{-\lambda t}f}$, then can study existence and uniqueness of the solution of the following problem:\\

\n $\left\{\begin{array}[c]{l}\displaystyle{v_{t} + v_{xxxx} + v_{xxx} + Rv_{xx} + \frac{1}{2}e^{\lambda t}v_{x}^{2} = \tilde{f};\quad (t, x) \in ]0, \infty[\times \mathbb{R}}\\
\quad\\
\displaystyle{v(t, x+ 2\pi) = v(t, x)}\\
\quad\\
\displaystyle{v(0, x) = v_{0}}\\
\end{array}\right.$ $\hfill { } {\color{red}(\mathbb{P})}$ \\
\quad\\

\n Let $\displaystyle{\mathfrak{H} := \mathbb{L}^{2}(0, T) = \{f : (0, T) \longrightarrow \mathbb{R}; \int_{0}^{T}|f(t)|^{2}dt < \infty, T > 0\}}$ and $H^{s}(0, T)$ be the classical Sobolev space of order $s$.\\

\n Let $\mathbb{X}$ be a Banach space with norm $\displaystyle{||.||_{\mathbb{X}}}$ and\\

\n - $\displaystyle{\mathbb{L}^{p}(0,T; \mathbb{X}) = \{f : (0, T) \longrightarrow \mathbb{X}; \int_{0}^{T}||f(t)||_{\mathbb{X}}^{p}dt < \infty ; p \in \mathbb{N}\}}$\\

\n - $\displaystyle{\mathbb{L}^{\infty}(0,T; \mathbb{X}) = \{f : (0, T) \longrightarrow \mathbb{X};\quad  _{_{_{_{t \in ]0, T[}}}}\!\!\!\!\!\!\!\!\!\!\!\!Sup.ess||f(t)||_{\mathbb{X}} < \infty\}}$\\

\n The main result of these results is the following theorem:\\

\begin{theorem}
\n Let $\displaystyle{\tilde{f} \in \mathbb{L}^{2}(0,T; \mathbb{H}^{1})}$ and $\displaystyle{v_{0} \in \mathbb{H}^{1}}$\\

\n Then there exist an uniqueness solution of the problem $\mathbb{P}$ satisfying\\

\n i) $\displaystyle{v \in \mathbb{L}^{2}(0,T; \mathbb{H}^{3})\cap\mathbb{L}^{\infty}(0,T; \mathbb{H}^{2})}$\\

\n ii) $\displaystyle{v_{t} \in \mathbb{L}^{2}(\Omega)}$ where $\Omega = ]0, T[\times]0, 2\pi[$\\

\end{theorem}

\n The proof of this theorem consists of serval lemmas. We begin by\\

\begin{lemma}  (a priori inequalities of $\mathbb{L}^{2}$ type for solutions of KdV-KS equation)\\

\n Let $\displaystyle{\tilde{f} \in \mathbb{L}^{2}(0,T; \mathbb{H}^{1})}$\\

\n i) If $\displaystyle{v_{0} \in \mathbb{H}^{1}}$ and $v$ is solution of $\mathbb{P}$ then we have\\

\n $\alpha$) $\displaystyle{||v_{x}(t)||^{2} + \int_{0}^{t}||v_{xxx}(s)||^{2}ds + (2\lambda - R^{2} -1) \int_{0}^{t}||v_{x}(s)||^{2}ds \leq C(v_{0}, \tilde{f})\quad \forall t \in ]0, T[}$\\

\n $\beta$) $\displaystyle{||v(t)||^{2} + 2\int_{0}^{t}||v_{xx}(s)||^{2}ds + (2\lambda - \frac{3}{4}) \int_{0}^{t}||v(s)||^{2}ds \leq C(v_{0}, \tilde{f})\quad \forall t \in ]0, T[}$\\

\n where $C(v_{0}, \tilde{f})$ is constant.\\

\n ii) If $\displaystyle{v_{0} \in \mathbb{H}^{2}}$ and $v$ is solution of $\mathbb{P}$ then we have\\

\n $\displaystyle{\int_{0}^{t}||v_{t}(s)||^{2}ds + ||v_{xx}(t)||^{2} + \lambda||v(t)||^{2} \leq C(v_{0}, \tilde{f})\quad \forall t \in ]0, T[}$\\

\n where $C(v_{0}, \tilde{f})$ is constant.\\

\n iii) If $\displaystyle{v_{0} \in \mathbb{H}^{3}}$ and $v$ is solution of $\mathbb{P}$ then we have\\

\n $\displaystyle{\int_{0}^{t}||v_{tx}(s)||^{2}ds + ||v_{xxx}(t)||^{2} + \lambda||v_{x}(t)||^{2} \leq C(v_{0}, \tilde{f})\quad \forall t \in ]0, T[}$\\

\n where $C(v_{0}, \tilde{f})$ is constant.\\
\end{lemma}
\n {\bf{\color{red} Proof}}\\

\n We work within the framework of the Hilbert space $\mathcal{H} = \mathbb{L}^{2}(0, 2\pi)$.\\

\n  The inner product and norm are given respectively by,\\

\n $\displaystyle{<u, v> = \int_{0}^{2\pi}u(x)v(x)dx}$ $\hfill { } {\color{red}(\mathbb{F}_{3})}$ \\

\n $\displaystyle{||u|| = \sqrt{\int_{0}^{2\pi}|u(x)|^{2}dx}}$ $\hfill { } {\color{red}(\mathbb{F}_{4})}$ \\

\n i) let $\displaystyle{v_{0} \in \mathbb{H}^{1}}$ and $v$ be solution of $\mathbb{P}$, to proof a priori inequality $\alpha$) multiplying ($\mathbb{F}_{2}$) by $-v_{xx}$ and integrating on $[0, 2\pi]$, then it easy to see from the periodic conditions that\\

\n $\displaystyle{\int_{0}^{2\pi}v_{xxx}v_{xx}dx = \frac{1}{2}\int_{0}^{2\pi}[v_{xx}^{2}]_{x}dx = \frac{1}{2}[v_{xx}^{2}]_{0}^{2\pi} = 0}$ $\hfill { } {\color{red}(\mathbb{F}_{5})}$ \\

\n and\\

$\displaystyle{\int_{0}^{2\pi}v_{x}^{2}v_{xx}dx = -\frac{1}{3}\int_{0}^{2\pi}[v_{x}^{3}]_{x}dx = -\frac{1}{3}[v_{x}^{3}]_{0}^{2\pi} = 0}$ $\hfill { } {\color{red}(\mathbb{F}_{6})}$ \\

\n Now by using Green formula and periodic conditions, we get\\

\n $\displaystyle{\frac{1}{2}\frac{d}{dt}||v_{x}||^{2} + ||v_{xxx}||^{2} + \lambda||v_{x}||^{2} = -R\int_{0}^{2\pi}v_{x}v_{xx}dx + \int_{0}^{2\pi}\tilde{f}_{x}v_{x}dx}$ $\hfill { } {\color{red}(\mathbb{F}_{7})}$ \\

\n we use the Cauchy-Schwarz inequality on the terms of second member of ($\mathbb{F}_{7}$) to obtain\\

\n $\displaystyle{-R\int_{0}^{2\pi}v_{x}v_{xx}dx| \leq \frac{1}{2}||v_{xxx}||^{2} + \frac{1}{2}R^{2}||v_{x}||^{2}}$\\

\n and\\

\n $\displaystyle{|\int_{0}^{2\pi}\tilde{f}_{x}v_{x}dx| \leq \frac{1}{2}||v_{x}||^{2} + \frac{1}{2}||\tilde{f}_{x}||^{2}}$\\

\n This implies the following a priori inequality\\

\n $\displaystyle{\frac{d}{dt}||v_{x}||^{2} + ||v_{xxx}||^{2} + (2\lambda - R^{2} -1)||v_{x}||^{2} \leq ||\tilde{f}_{x}||^{2}}$ $\hfill { } {\color{red}(\mathbb{F}_{8})}$ \\

\n integrating ($\mathbb{F}_{8}$) of $0$ to $t$ to get\\

\n $\displaystyle{||v_{x}(t)||^{2} + \int_{0}^{t}||v_{xxx}(s)||^{2}ds + (2\lambda - R^{2} -1)\int_{0}^{t}||v_{x}(s)||^{2}ds \leq C(v_{0}, \tilde{f})}$ $\hfill { } {\color{red}(\mathbb{F}_{9})}$ \\

\n This finish the proof of $\alpha$) of i) of this lemma.\\

\begin{remark}
\quad\\

\n If we choose $\lambda \geq \frac{R^{2}+1}{2}$ we have\\

\n a) $\displaystyle{v_{x} \in \mathbb{L}^{\infty}(0, T; \mathbb{L}^{2})}$\\

\n b) $\displaystyle{v_{xxx} \in \mathbb{L}^{2}(0, T; \mathbb{L}^{2})}$\\

\n c) $\displaystyle{v_{t} \in \mathbb{L}^{2}(0, T; \mathbb{L}^{2})}$\\

\n d) $\displaystyle{v_{xx} \in \mathbb{L}^{2}(0, T; \mathbb{L}^{2})}$\\

\n e) $\displaystyle{v_{x} \in \mathbb{L}^{2}(0, T; \mathbb{H}^{2})}$ and $\displaystyle{v_{x} \in \mathbb{L}^{2}(0, T; \mathbb{L}^{\infty})}$\\
\end{remark}

\n To proof a priori inequality $\beta$), multiplying ($\mathbb{F}_{2}$) by $v$ and integrating on $[0, 2\pi]$, then it easy to see from the periodic conditions that\\

\n $\displaystyle{\int_{0}^{2\pi}v_{xxx}vdx = [v_{xx}v]_{0}^{2\pi} -\int_{0}^{2\pi}v_{xx}v_{x}dx = [v_{xx}v]_{0}^{2\pi} - \frac{1}{2}[v_{x}^{2}]_{0}^{2\pi}= 0}$ $\hfill { } {\color{red}(\mathbb{F}_{9})}$ \\

\n Now by using Green formula and once more periodic conditions, we get\\

\n $\displaystyle{\frac{1}{2}\frac{d}{dt}||v||^{2} + ||v_{xx}||^{2} + \lambda||v||^{2} = R||v_{x}||^{2} -\frac{1}{2}e^{\lambda t}\int_{0}^{2\pi}v_{x}^{2}vdx + \int_{0}^{2\pi}\tilde{f}_{x}vdx}$ $\hfill { } {\color{red}(\mathbb{F}_{10})}$ \\

\n As\\

\n $\displaystyle{|e^{\lambda t}\int_{0}^{2\pi}v_{x}^{2}vdx| \leq e^{\lambda t}||v_{x}||_{\infty}||v_{x}||||v|| \leq \frac{1}{2}e^{2\lambda t}||v_{x}||_{\infty}^{2}||v_{x}||^{2} + \frac{1}{2}||v||^{2}}$\\

\n and\\

\n $\displaystyle{|\int_{0}^{2\pi}\tilde{f}vdx| \leq \frac{1}{2}||v||^{2} + \frac{1}{2}||\tilde{f}||^{2}}$\\

\n we deduce the following a priori inequality\\

\n $\displaystyle{\frac{d}{dt}||v||^{2} + 2||v_{xx}||^{2} + (2\lambda - \frac{3}{2})||v||^{2} \leq ||v_{x}||^{2}[2R + \frac{1}{2}||v_{x}||_{\infty}^{2}] + ||\tilde{f}_{x}||^{2}}$ $\hfill { } {\color{red}(\mathbb{F}_{11})}$ \\

\n integrating ($\mathbb{F}_{11}$) of $0$ to $t$ and as (see above remark) $\displaystyle{v_{x} \in \mathbb{L}^{2}(0, T; \mathbb{L}^{\infty})\cap \mathbb{L}^{\infty}(0, T; \mathbb{L}^{2})}$ we deduce that\\

\n $\displaystyle{||v(t)||^{2} + 2\int_{0}^{t}||v_{xx}(s)||^{2}ds + (2\lambda - \frac{3}{2})\int_{0}^{t}||v(s)||^{2}ds \leq C(v_{0}, \tilde{f})}$ $\hfill { } {\color{red}(\mathbb{F}_{12})}$ \\

\n This finish the proof of the property $\beta$).\\

\begin{remark}
\quad\\

\n If we choose $\lambda \geq \frac{3}{4}$ we deduce that\\

\n a) $\displaystyle{ u \in \mathbb{L}^{\infty}(0, T; \mathbb{L}^{2})}$\\

\n b) $\displaystyle{ u \in \mathbb{L}^{2}(0, T; \mathbb{L}^{2})}$\\
\end{remark}

\n ii) Let $\displaystyle{v_{0} \in \mathbb{H}^{2}}$ and $v$ be solution of $\mathbb{P}$, to proof a priori inequality ii), multiplying ($\mathbb{F}_{2}$) by $v_{t}$ and integrating on $[0, 2\pi]$ and by using Green formula with once more periodic conditions, we get\\

\n $\displaystyle{||v_{t}||^{2} + \frac{1}{2}\frac{d}{dt}||v_{xx}||^{2} + \frac{\lambda}{2}\frac{d}{dt}||v||^{2} = -\int_{0}^{2\pi}v_{xxx}v_{t}dx - R\int_{0}^{2\pi}v_{xx}v_{t}dx -}$\\
$\displaystyle{\frac{1}{2}e^{\lambda t}\int_{0}^{2\pi}v_{x}^{2}v_{t}dx +
\int_{0}^{2\pi}\tilde{f}_{x}v_{t}dx}$ $\hfill { } {\color{red}(\mathbb{F}_{13})}$ \\

\n Now by using the following inequalities\\

\n 1) $\displaystyle{|\int_{0}^{2\pi}v_{xxx}v_{t}dx| \leq \frac{1}{8}||v_{t}||^{2} + 2||v_{xxx}||^{2}}$\\

\n 2) $\displaystyle{|R\int_{0}^{2\pi}v_{xx}v_{t}dx| \leq \frac{1}{8}||v_{t}||^{2} + 2R^{2}||v_{xx}||^{2}}$\\

\n 3) $\displaystyle{|e^{\lambda t}\int_{0}^{2\pi}v_{x}^{2}v_{t}dx| \leq e^{\lambda t}||v_{x}||_{\infty}||v_{x}||||v_{t}||\leq e^{2\lambda t}||v_{x}||_{\infty}^{2}||v_{x}||^{2} +
\frac{1}{4}||v_{t}||^{2}}$\\

\n 4) $\displaystyle{\int_{0}^{2\pi}\tilde{f}_{x}v_{t}dx| \leq \frac{1}{8}||v_{t}||^{2} + 2||\tilde{f}||^{2}}$\\

\n we get the following inequality\\

\n $\displaystyle{||v_{t}||^{2} + \frac{d}{dt}||v_{xx}||^{2} + \lambda \frac{d}{dt}||v||^{2} \leq 4||v_{xxx}||^{2} + 4R^{2}||v_{xx}||^{2} + }$\\

\n $\displaystyle{e^{\lambda t}||v_{x}||_{\infty}^{2}||v_{x}||^{2} + 4||\tilde{f}_{x}||^{2}}$ $\hfill { } {\color{red}(\mathbb{F}_{14})}$ \\

\n integrating ($\mathbb{F}_{14}$) of $0$ to $t$ \\

\n As $\displaystyle{u_{xx} \in \mathbb{L}^{2}(0, T; \mathbb{L}^{2}),u_{xxx} \in \mathbb{L}^{2}(0, T; \mathbb{L}^{2})}$ and $\displaystyle{ v_{x} \in \mathbb{L}^{2}(0, T; \mathbb{L}^{\infty})\cap \mathbb{L}^{\infty}(0, T; \mathbb{L}^{2})}$ we deduce that\\

\n $\displaystyle{\int_{0}^{t}||v_{t}(s)||^{2}ds + ||v_{xx}(t)||^{2} + \lambda||v(t)||^{2} \leq C(v_{0}, \tilde{f})}$ $\hfill { } {\color{red}(\mathbb{F}_{15})}$ \\

\n This finish the proof of the property ii).\\

\begin{remark}
\quad\\

\n If we choose $\lambda \geq 0$ we deduce that\\

\n a) $\displaystyle{ u_{t} \in \mathbb{L}^{2}(0, T; \mathbb{L}^{2})}$\\

\n b) $\displaystyle{ u_{xx} \in \mathbb{L}^{\infty}(0, T; \mathbb{L}^{2})}$\\
\end{remark}

\n iii) let $\displaystyle{v_{0} \in \mathbb{H}^{3}}$ and $v$ be solution of $\mathbb{P}$, to proof a priori inequality ii), multiplying ($\mathbb{F}_{2}$) by $v_{txx}$ and integrating on $[0, 2\pi]$ and by using Green formula with once more periodic conditions, we get\\

\n $\displaystyle{||v_{tx}||^{2} + \frac{1}{2}\frac{d}{dt}||v_{xxx}||^{2} + \frac{\lambda}{2}\frac{d}{dt}||v_{x}||^{2} = -\int_{0}^{2\pi}v_{xxxx}v_{tx}dx - R\int_{0}^{2\pi}v_{xxx}v_{tx}dx - }$\\
$\displaystyle{e^{\lambda t}\int_{0}^{2\pi}v_{x}v_{xx}v_{tx}dx + \int_{0}^{2\pi}\tilde{f}_{x}v_{tx}dx}$ $\hfill { } {\color{red}(\mathbb{F}_{16})}$ \\

\n Now by using the following inequalities\\

\n 1) $\displaystyle{|\int_{0}^{2\pi}v_{xxxx}v_{tx}dx| \leq \frac{1}{8}||v_{tx}||^{2} + 2||v_{xxxx}||^{2}}$\\

\n 2) $\displaystyle{|R\int_{0}^{2\pi}v_{xxx}v_{tx}dx| \leq \frac{1}{8}||v_{tx}||^{2} + 2R^{2}||v_{xxx}||^{2}}$\\

\n 3) $\displaystyle{|e^{\lambda t}\int_{0}^{2\pi}v_{x}v_{xx}v_{tx}dx| \leq e^{\lambda t}||v_{x}||_{\infty}||v_{xx}||||v_{tx}||\leq 2e^{2\lambda t}||v_{x}||_{\infty}^{2}||v_{xx}||^{2} +\frac{1}{8}||v_{tx}||^{2}}$\\

\n 4) $\displaystyle{\int_{0}^{2\pi}\tilde{f}_{x}v_{tx}dx| \leq \frac{1}{8}||v_{tx}||^{2} + 2||\tilde{f}_{x}||^{2}}$\\

\n we get the following inequality\\

\n $\displaystyle{||v_{tx}||^{2} + \frac{d}{dt}||v_{xxx}||^{2} + \lambda \frac{d}{dt}||v_{x}||^{2} \leq 4||v_{xxxx}||^{2} + 4R^{2}||v_{xxx}||^{2} +}$\\
$\displaystyle{ 4e^{\lambda t}||v_{x}||_{\infty}^{2}||v_{xx}||^{2} + 4||\tilde{f}_{x}||^{2}}$ $\hfill { } {\color{red}(\mathbb{F}_{17})}$ \\

\n Now, we observe that\\

\n $\displaystyle{v_{xxxx} = -v_{t} - v_{xxx} - Rv_{xx} -\lambda v - \frac{1}{2}e^{\lambda t}v_{x}^{2} + \tilde{f}}$\\

\n we deduce\\

\n $\displaystyle{||v_{xxxx}||^{2}\leq 6[||v_{t}||^{2}+ ||v_{xxx}||^{2}+R^{2}||v_{xx}||^{2}+ \lambda^{2}||v||^{2}] + }$\\
$\displaystyle{6[\frac{1}{4}e^{2\lambda t}||v_{x}^{2}||^{2} + ||\tilde{f}||^{2}]}$ $\hfill { } {\color{red}(\mathbb{F}_{18})}$ \\

\n $\displaystyle{||v_{x}^{2}||^{2} \leq ||v_{x}||_{\infty}^{2}||v_{x}||^{2}}$ $\hfill { } {\color{red}(\mathbb{F}_{19})}$ \\

\n integrating ($\mathbb{F}_{17}$) of $0$ to $t$ then we obtain\\

\n $\displaystyle{\int_{0}^{t}||v_{tx}(s)||^{2}ds + ||v_{xxx}(t)||^{2} + \lambda||v_{x}(t)||^{2} \leq ||v_{0_{_{xxx}}}||^{2} + \lambda||v_{0_{_{x}}}||^{2} + 4\int_{0}^{t}||v_{xxxx}(s)||^{2}ds +}$\\
$\displaystyle{4R^{2}\int_{0}^{t}||v_{xxx}(s)||^{2}ds + 4\int_{0}^{t}||\tilde{f}_{x}(s)||^{2}ds + 4\int_{0}^{t}e^{2\lambda s}||v_{x}(s)||_{\infty}^{2}||v_{xx}(s)||^{2}ds}$ $\hfill { } {\color{red}(\mathbb{F}_{20})}$ \\

\n From the inequalities established in i) and ii) and the asymptions on $v_{0}$ and $f$ we know that all terms of second member of ($\mathbb{F}_{20}$) are bounded excepted the terms $\displaystyle{\int_{0}^{t}||v_{xxxx}(s)||^{2}ds}$ and $\displaystyle{\int_{0}^{t}e^{2\lambda s}||v_{x}(s)||_{\infty}^{2}||v_{xx}(s)||^{2}ds}$ they are bounded by:\\

\n $\displaystyle{\int_{0}^{t}e^{2\lambda s}||v_{x}(s)||_{\infty}^{2}||v_{xx}(s)||^{2}ds \leq e^{2\lambda t}||v_{xx}||_{\mathbb{L}^{\infty}(0, T; \mathbb{L}^{2})}||v_{x}||_{\mathbb{L}^{2}(0, T; \mathbb{L}^{\infty})}}$\\

\n and\\

\n $\displaystyle{\int_{0}^{t}||v_{xxxx}(s)||^{2}ds \leq 6[\int_{0}^{t}||v_{t}(s)||^{2}ds + \int_{0}^{t}||v_{xxx}(s)||^{2}ds + R^{2}\int_{0}^{t}||v_{xx}(s)||^{2}ds] + }$\\
$\displaystyle{6[\lambda^{2}\int_{0}^{t}||v(s)||^{2}ds] +\frac{1}{4}\int_{0}^{t}e^{2\lambda s}||v_{x}^{2}(s)||^{2}ds + \int_{0}^{t}||\tilde{f}(s)||^{2}ds]}$\\

\n As $v_{x}$ is bounded in $\mathbb{L}^{2}(0, T; \mathbb{L}^{\infty})$ and $v_{xx}$ is bounded in $\mathbb{L}^{\infty}(0, T; \mathbb{L}^{2})$ we deduce that the term $\displaystyle{\int_{0}^{t}e^{2\lambda s}||v_{x}(s)||_{\infty}^{2}||v_{xx}(s)||^{2}ds}$ is bounded for all $t \in ]0, T[$\\

\n Now from ($F_{19}$), i) and ii) we deduce also that the term $\displaystyle{\int_{0}^{t}||v_{xxxx}(s)||^{2}ds}$ is bounded for all $t \in ]0, T[$\\

\n Then\\

\n $\displaystyle{\int_{0}^{t}||v_{tx}(s)||^{2}ds + ||v_{xxx}(t)||^{2} + \lambda||v_{x}(t)||^{2} \leq C(v_{0}, \tilde{f})}$ \quad $\forall t \in ]0, T[$ $\hfill { } {\color{red}(\mathbb{F}_{21})}$ \\

\n This finish the proof of the property iii).\\

\begin{lemma}
\quad\\

\n If $v$ is solution of the problem $\mathbb{P}$, then we have uniqueness of the solution.\\
\end{lemma}

\n {\bf{\color{red} Proof}}\\

\n We consider two solutions $u$ and $v$ of the problem $\mathbb{P}$ and we put $u = v - w$ then $u$ is solution of the following problem\\

\n $\left\{\begin{array}[c]{l}\displaystyle{u_{t} + u_{xxxx} + u_{xxx} + Ru_{xx} + \frac{1}{2}e^{\lambda t}u_{x}^{2} = 0;\quad (t, x) \in ]0, \infty[\times \mathbb{R}}\\
\quad\\
\displaystyle{u(t, x + 2\pi) = u(t, x)}\\
\quad\\
\displaystyle{u(0, x) = 0 ; \quad \forall x \in \mathbb{R}}\\
\end{array}\right.$ $\hfill { } {\color{red}(\mathbb{P}_{0})}$ \\
\quad\\

\n Multiplying equation of the problem ($\mathbb{P}_{0}$) in $L^{2}(0, 2\pi)$ by $-u_{xx}$, by following the procedure used in the above lemma we obtain\\

\n $\displaystyle{\frac{1}{2} \frac{d}{dt}||u_{x}||^{2} + ||u_{xxx}||^{2} + \lambda||u_{x}||^{2} = -R\int_{0}^{2\pi}u_{x}u_{xxx}dx \quad + }$\\

\n $\displaystyle{\frac{1}{2}e^{\lambda t}\int_{0}^{2\pi}u_{x}(v_{x} + w_{x})u_{xx}dx}$$\hfill { } {\color{red}(\mathbb{F}_{22})}$ \\

\n Now we use the following inequalities\\

\n 1) $\displaystyle{e^{\lambda t}\int_{0}^{2\pi}u_{x}(v_{x} + w_{x})u_{xx}dx \leq \frac{1}{2}e^{2\lambda t}||v_{x} + w_{x}||_{\infty}^{2}||u_{x}||^{2} + \frac{1}{2}||u_{xx}||^{2}}$\\

\n 2)$\displaystyle{|R\int_{0}^{2\pi}u_{x}u_{xxx}dx|  \leq \frac{1}{2}||u_{xxx}||^{2} + \frac{1}{2}R^{2}||u_{x}||^{2}}$\\

\n to deduce\\

\n $\displaystyle{\frac{d}{dt}||u_{x}||^{2} + ||u_{xxx}||^{2} - \frac{1}{2}||u_{xx}||^{2} + (2\lambda - R^{2})||u_{x}||^{2} \leq }$\\

\n $\displaystyle{\frac{1}{2}e^{2\lambda t}||v_{x} + w_{x}||_{\infty}^{2}||u_{x}||^{2} }$ $\hfill { } {\color{red}(\mathbb{F}_{23})}$ \\

\n Multiplying equation of problem ($\mathbb{P}_{0}$) in $L^{2}(0, 2\pi)$ by $u$ then we get\\

\n $\displaystyle{\frac{d}{dt}||u||^{2} + ||u_{xx}||^{2} + \lambda||u||^{2} =  - R\int_{0}^{2\pi}u_{xx}u_{x}dx - \frac{1}{2}e^{\lambda t}\int_{0}^{2\pi}u_{x}(v_{x} + w_{x})udx}$\\

\n As\\

\n $\displaystyle{e^{\lambda t}\int_{0}^{2\pi}u_{x}(v_{x} + w_{x})u_{xx}dx \leq \frac{1}{2}e^{2\lambda t}||v_{x} + w_{x}||_{\infty}^{2}||u_{x}||^{2} + \frac{1}{2}||u_{xx}||^{2}}$\\

\n and\\

\n $\displaystyle{|R\int_{0}^{2\pi}u_{xx}udx|  \leq \frac{1}{2}||u_{xx}||^{2} + \frac{1}{2}R^{2}||u||^{2}}$\\

\n then we have\\

\n $\displaystyle{\frac{d}{dt}||u||^{2} + ||u_{xx}||^{2}  + (2\lambda - 1 - R^{2})||u||^{2} \leq e^{2\lambda t}||v_{x} + w_{x}||_{\infty}^{2}||u_{x}||^{2} }$ $\hfill { } {\color{red}(\mathbb{F}_{24})}$ \\

\n By adding ($\mathbb{F}_{23}$) and ($\mathbb{F}_{24}$) we obtain\\

\n $\displaystyle{\frac{d}{dt}||u||_{H^{1}}^{2} + ||u_{xxx}||^{2}  + \frac{1}{2}||u_{xx}||^{2} + (2\lambda - R^{2})||u_{x}||^{2} + (2\lambda -1 - R^{2})||u||^{2} \leq }$\\

\n $\displaystyle{\frac{3}{2}e^{2\lambda t}||v_{x} + w_{x}||_{\infty}^{2}||u_{x}||^{2} }$ $\hfill { } {\color{red}(\mathbb{F}_{25})}$ \\

\n by taking $\displaystyle{\lambda \geq \frac{R^{2}+1}{2}}$, we deduce that\\

\n $\displaystyle{\frac{d}{dt}||u||_{H^{1}}^{2} \leq \frac{3}{2}e^{2\lambda t}||v_{x} + w_{x}||_{\infty}^{2}||u_{x}||^{2} }$\\

\n and consequently \\

\n $\displaystyle{\frac{d}{dt}||u||_{H^{1}}^{2} \leq \frac{3}{2}e^{2\lambda t}||v_{x} + w_{x}||_{\infty}^{2}||u||_{H^{1}}^{2} }$$\hfill { } {\color{red}(\mathbb{F}_{26})}$ \\

\n Now, we use the Gronwall lemma to deduce that $\displaystyle{||u||_{H^{1}}^{2} = 0}$ and this implies that $u = 0$, hence we have uniqueness of the solution of the problem $\mathbb{P}$\\

\n Now, we use Galerkin procedure for the proof of the existence of solution of KdV-KS equation based on its variational formulation\\

\n Let $\displaystyle{\mathcal{V}_{m} = span\{u_{1}, u_{2}, ..., u_{m}\}}$ be space spanned by the $m$ first eigenvectors  of the operator $\Delta$ defined by $\Delta v = - v_{xx}$ with domain:\\

\n  $\displaystyle{D(\Delta) = \{ v \in \mathbb{H}_{loc}^{2}; v(x + 2\pi) = v(x) \quad \forall x \in \mathbb{R}\}}$\\

\n Defining an approximation of the solution of problem $\mathbb{P}$ by\\

\n $\displaystyle{v^{m} = \sum_{j=1}^{m}\alpha_{j}^{m}u_{j}}$ where $\alpha_{j}^{m}$ satisfy the conditions of the following system\\

\n $\left\{\begin{array}[c]{l}\displaystyle{<v_{t}^{m}, u_{j}> + <v_{xxxx}^{m}, u_{j}> + <v_{xxx}^{m}, u_{j}> + R<v_{xx}^{m}, u_{j}> + \lambda<v^{m}, u_{j}> +}\\
\quad\\
\displaystyle{\frac{1}{2}e^{\lambda t}<(v_{x}^{m})^{2}, u_{j}> = <\tilde{f}, u_{j}>; j=1, ..., m}\\
\quad\\
\displaystyle{v^{m}(0) = v_{0}^{m} \in \mathbb{H}^{2}(0, 2\pi)\quad\quad\quad\quad \quad\quad \quad\quad\quad\quad \quad\quad \quad\quad\quad\quad \quad\quad \quad\quad {\color{red} (\mathbb{F}_{27})}}\\
\quad\\
\displaystyle{v_{0}^{m} \longrightarrow v_{0}\quad in \quad \mathbb{H}^{1}(0, 2\pi)}\\
\quad\\
\displaystyle{\exists c > 0; ||v_{0}^{m}||_{\mathbb{H}^{2}} \leq c||v_{0}||_{\mathbb{H}^{1}}}\\
\end{array}\right.$ \\

\n As\\

 \n $\displaystyle{<v_{xxx}^{m}, u_{j}> = <\sum_{j=1}^{m}\alpha_{j}^{m}(t)\frac{d^{3}}{dx^{3}}u_{j}, u_{j}> = \sum_{k=1}^{m}\alpha_{j}^{m}(t)\int_{0}^{2\pi}\frac{d^{3}}{dx^{3}}(u_{k})u_{j}dx}$\\

\n  $\displaystyle{-\sum_{k=1}^{m}\lambda_{k}\alpha_{j}^{m}(t)\int_{0}^{2\pi}\frac{d}{dx}(u_{k})u_{j}dx}$\\

 \n where $\lambda_{k}$ is an eigenvalue of operator $\Delta$ associated to eigenvectors $u_{k}$.\\

 \n Then if we put $\displaystyle{\gamma_{kj} = \int_{0}^{2\pi}\frac{d}{dx}u_{k}u_{j}dx}$, we can write\\

\n $\displaystyle{<v_{xxx}^{m}, u_{j}> = \sum_{k=1}^{m}-\lambda_{k}\gamma_{kj}\alpha_{k}^{m}(t)}$

\n Now, for nonlinear term, we can write \\

\n  $\displaystyle{<(v_{x}^{m})^{2}, u_{j}> = \int_{0}^{2\pi} [\sum_{i=1}^{m}\alpha_{i}^{m}(t)\frac{d}{dx}u_{i}]^{2}u_{j}dx = \int_{0}^{2\pi} [\sum_{i=1}^{m}\sum_{l=1}^{m}c_{il}\alpha_{i}^{m}(t)\alpha_{l}^{m}(t)\frac{d}{dx}u_{i}\frac{d}{dx}u_{l}]u_{j}dx }$\\
 $\displaystyle{= \sum_{i=1}^{m}\sum_{l=1}^{m}\gamma_{ilj}\alpha_{i}^{m}(t)\alpha_{l}^{m}(t)}$ where $\displaystyle{ \gamma_{ilj} = \frac{d}{dx}u_{i}\frac{d}{dx}u_{l}]u_{j}dx}$\\

\n The coefficients $\alpha_{i}^{m}$ satisfy the following system:\\

\n $\left\{\begin{array}[c]{l}\displaystyle{\frac{d}{dt}\alpha_{j}^{m} + (\lambda + \lambda_{j}^{2} - R\lambda_{j})\alpha_{j}^{m} - \sum_{k=1}^{m}\lambda_{k}\gamma_{kj}\alpha_{k}^{m}}\\
\quad\\
\displaystyle{+ \quad\frac{1}{2}e^{\lambda t}\sum_{i=1}^{m}\sum_{l=1}^{m}\gamma_{ilj}\alpha_{i}^{m}(t)\alpha_{l}^{m}(t) = \tilde{f}_{j}(t); j=1, ..., m} \quad\quad\quad\quad \quad\quad\quad{\color{red} (\mathbb{F}_{28})}\\
\quad\\
\displaystyle{\alpha_{j}^{m}(0) = < v_{0}^{m}, u_{j}> }\\
\end{array}\right.$ \\

\n where $\displaystyle{\tilde{f}_{j}(t) =  < \tilde{f}, u_{j}>}$\\

\n The non linear system ($\mathbb{F}_{28}$) is loc-Lipschitzian on $\mathcal{V}_{m}$ hence its solution exists and it is maximal with continuous derivative on $[0, t_{m}[$ ($0 \leq t_{m} < + \infty$) \\

 \n Under the two following conditions:\\

\n (i) $\displaystyle{\tilde{f} \in \mathbb{L}^{2}(0, T, \mathbb{H}^{1})}$\\

\n (ii) $\displaystyle{v_{0}^{m} \in \mathbb{H}^{2}; \exists c > 0 , ||v_{0}^{m} ||_{\mathbb{H}^{2}} \leq c|| v_{0} ||_{\mathbb{H}^{1}}}$\\

\n  we present some properties of $v^{m}$:\\

 \n 1) $\displaystyle{v^{m}, v_{x}^{m}, v_{t}^{m}}$ satisfy a priori inequalities of the first lemma of this annex and from (ii) we deduce that $\displaystyle{||v_{0}^{m}|| \leq C(v_{0}, \tilde{f})}$ uniformly and $t_{m} = T$\\

 \n 2) For $\displaystyle{\lambda \geq sup(\frac{R^{2}+1}{2}, \frac{3}{4})}$ and $m \longrightarrow +\infty$ we have\\

\n  $\alpha$) $v^{m}$ is in a bounded set of $\mathbb{L}^{\infty}(0, T, \mathbb{H}^{3})$ \\

 \n $\beta$) $v^{m}$ is in a bounded set of $\mathbb{L}^{2}(0, T, \mathbb{H}^{4})$ \\

 \n $\gamma$) $v_{t}^{m}$ is in a bounded set of $\mathbb{L}^{\infty}(0, T, \mathbb{H}^{3})$ \\

 \n $\delta$) $v_{tx}^{m}$ is in a bounded set of $\mathbb{L}^{\infty}(0, T, \mathbb{H}^{3})$ \\

.\n $\hfill { } {\color{red}(\clubsuit)}$ \\

\n These properties allow us to pass to the limit if we found one compactly property to pass to the limit in the nonlinear term.\\
\n Considering  the Sobolev space $\displaystyle{\mathbb{H}^{r,s}(Q)= \mathbb{L}^{2}(0, T; \mathbb{H}^{r})\cap \mathbb{H}^{s}(0, T; \mathbb{L}^{2})}$ where \\$Q = [0, T]\times [0, 2\pi]$. \\

\n From the above properties we deduce that $v_{x}^{m}$ is in a bounded set of $\mathbb{H}^{3,1}(Q)$ hence it is  in a bounded set of $\mathbb{H}^{1}(Q)$. Now by the regularity results of Sobolev that the injection of $\mathbb{H}^{1}(Q)$ in $\mathbb{L}^{4}(Q)$ is compact allow us to extract a weakly convergent subsequence which we continue to denote by $v^{m}$ such that\\

\n a) $v^{m} \longrightarrow v $ in $\mathbb{L}^{2}(0, T; \mathbb{H}^{4})$ weakly.\\

\n b) $v_{t}^{m} \longrightarrow v_{t} $ in $\mathbb{L}^{2}(Q)$ weakly.\\

\n c) $v^{m} \longrightarrow v $ in $\mathbb{L}^{\infty}(0, T; \mathbb{H}^{3})$ weakly*.\\

\n d) $v_{x}^{m} \longrightarrow v $ in $\mathbb{L}^{4}(Q)$ strongly.\\

\n Then we can to pass to limit in first equation of $\mathbb{F}_{27}$ to get\\

\n $\displaystyle{<v_{t}, u_{j}> + <v_{xxxx}, u_{j}> + <v_{xxx}, u_{j}> + R<v_{xx}, u_{j}> + \lambda <v, u_{j}> }$\\

\n $\displaystyle{+\frac{1}{2}e^{\lambda t}<v_{x}^{2}, u_{j}> = <\tilde{f}, u_{j}>, \quad \forall j}$ $\hfill {} {\color{red}(\mathbb{F}_{29})}$\\

\n and since the system of eigenvectors of $\Delta$ is dense in $\mathbb{H}^{3}$ we deduce that\\

\n $\displaystyle{<v_{t}, u> + <v_{xxxx}, u> + <v_{xxx}, u> + R<v_{xx}, u> + \lambda <v, u> }$\\

\n $\displaystyle{+\frac{1}{2}e^{\lambda t}<v_{x}^{2}, u> = <\tilde{f}, u>, \quad \forall u \in \mathbb{H}^{3}}$ $\hfill {} {\color{red}(\mathbb{F}_{30})}$\\

\n As $v$ is continuous of $[0, T]$ into $\mathbb{H}^{1}$ and from the above property c) with  the fact that $v_{0}^{m}$ is uniformly bounded, we deduce that $v(0, x) = v_{0} \in \mathbb{H}^{1}$  thus finish the proof of existence of solution of KdV-KS equation.\\

\begin{remark}
\quad\\

\n 1) To pass to the limit in the nonlinear term, we can observed (as $v_{t}^{m}$ is in bounded set of $\mathbb{L}^{2}(0, T; \mathbb{L}^{2})$) that $v^{m}$ is in bounded set of $\mathbb{H}^{0,1}(Q)$ hence $v_{x}^{m}$ is in bounded set of $\mathbb{H}^{-1,1}(Q)$. Let $\mathbb{D} = \{(x, y); y = ax + b, a \in \mathbb{R}, b \in \mathbb{R}\}$ be the line such that $(3, 0) \in \mathbb{D}$ and $(-1, 1) \in \mathbb{D}$. The point $(x, x) \in \mathbb{D}$ is  $(\frac{3}{5}, \frac{3}{5})$, from this geometric interpretation we deduce that the space of interpolation between $\mathbb{H}^{3,0}(Q)$ and $\mathbb{H}^{-1,1}(Q)$ is $\mathbb{H}^{\frac{3}{5},\frac{3}{5}}(Q)$ hence $v_{x}^{m}$ is in bounded set of $\mathbb{H}^{\frac{3}{5},\frac{3}{5}}(Q)$. Now as the injection of  $\mathbb{H}^{\frac{3}{5},\frac{3}{5}}(Q)$ in $\mathbb{H}^{\frac{1}{2},\frac{1}{2}}(Q)$ is compact and the injection of $\mathbb{H}^{\frac{1}{2},\frac{1}{2}}(Q)$ in $\mathbb{L}^{4}(Q)$ is continuous then the injection of $\mathbb{H}^{\frac{3}{5},\frac{3}{5}}(Q)$ in $\mathbb{L}^{4}(Q)$ is also compact and this complete our discussion on existence of solutions of KdV-KS equation.\\

\n 2) An other generalisation of KS equation is given by:\\

\n $\left \{ \begin{array}[c]{l} \displaystyle{u_{t } + \nu u_{xxxx } + \delta u_{xxx } + u_{xx} + uu_{x} = 0}\\
\quad\\
u(t, x) = u(t, x + 2\pi)\\
\quad\\
u(0, x) = u_{0}(x)\\
\end{array}\right .$ $\hfill {} {\color{red}(\mathbb{F}_{31})}$\\
\n where $\displaystyle{\nu = ( \frac{2\pi}{L})^{2}}$ is a positive parameter that decreases as the system size $L$ increases and the parameter $\delta$ measures dispersive effects.\\ 
\n It is well known that the equation {\color{red}($\mathbb{F}_{31}$)} is of the active-dissipative type and instabilities are present depending on the value of $\nu$ and if $\nu > 1$ the zero solution, representing a flat film, is unique. However, when $\nu < 1$ the zero solution is linearly unstable and bifurcates into nonlinear states including steady states, travelling waves and solutions exhibiting spatiotemporal chaos. Some of these solutions are stable, and others are unstable. It is established that sufficiently large values of $\delta$ act to regularise the dynamics (even chaotic ones) into nonlinear travelling wave pulses - see Kawahara in {\color{blue}[[15]]} and Kawahara-Toh in {\color{blue}[[16]]},  {\color{blue}[[17]]}. However, in a regime of moderate values of $\delta$ travelling waves or pulses appear to be randomly interacting with each other giving rise to what is widely known as weak/dissipative turbulence (in the ``Manneville sense"  {\color{blue}[[26]]}) for a weak interaction theory between pulses that are sufficiently separated. \\
\end{remark}
\n \begin{center}
{\color{blue}{\bf Time evolution of the equation {\color{red}($\mathbb{F}_{31}$)} for  $\nu = 0.01 \,\, (L = 20\pi)$}}\\
\end{center}
\n \includegraphics[scale=0.62]{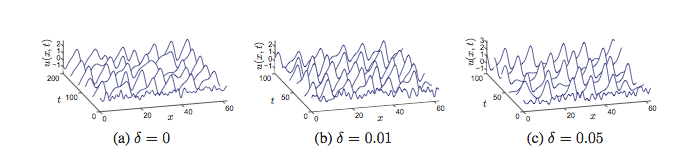}\\ 
\n {\color{red}(a)} shows the chaotic behaviour in the absence of dispersion\\
\n {\color{red} (b)} shows the weak/dissipative turbulent behaviour for small values of $\delta$\\
\n {\color{red}(c)}  shows the chaotic regularisation with relatively large values of dispersion.\\

\section{{\color{red}{\bf A study of the nonhomogeneous boundary value problem for KdV-KS equation in quarter plane: (2016)}}} 

\n In 2016, Jing Li,  Bing-Yu Zhang and Zhixiong Zhang  in {\color{blue} [[25]]} was studied the local well-posedness and the global well-posedness of the non homogeneous initial
boundary value problem of KdV-KS equation in quarter plane.\\
\n They was considered the following problem:\\
\begin{equation}
\left \{ \begin{array} {c} \displaystyle{u_t + u_{xxx} + \delta u_{xxx} + u_{xx} + uu_{x} = 0, \quad (t, x) \in (0, T)\times \mathbb{R}^{+}}\\
\quad\\
u(0, x) = \phi(x) \quad\quad\quad \quad\quad\quad\quad \quad\quad \quad \quad x\in \mathbb{R}^{+}\quad  \quad\quad\quad\\
\quad\\
u(t, 0) = h_{1}(t) \quad u_{x}(t, 0) = h_{2}(t) \quad \quad \quad\quad t \in (0, T),\quad\quad\\
\end{array} \right. \quad
\end{equation}

\n where $T > 0$, $\displaystyle{\delta \in \mathbb{R}}$ and $\displaystyle{\phi \in H^{s}(\mathbb{R}^{+})}$.\\

\n They was conjectured that the well-posedness will fail when $s < -2.$.\\

\n Some their main results was based on $s$-compatibility conditions.\\

\begin{definition}

\n Let $T > 0$ and $s  0$ be given.\\

\n $\displaystyle{(\phi, h_{1}, h_{2}) \in H^{s}(\mathbb{R}^{+}) \times H^{\frac{s}{4} + \frac{3}{8}}(0, T) \times H^{\frac{s}{4} +\frac{1}{8} }(0, T)}$ is said to be $s$-compatible if \\

\n (i) When $\displaystyle{s - 4[\frac{s}{4}] \leq \frac{1}{2}}$, then\\
\n $\displaystyle{\phi_{k}(0) = h_{1}^{(k)}(0), \quad \phi_{k}^{'}(0) = h_{2}^{(k)}(0), \quad\quad k = 0, 1, ......, [\frac{s}{4}]-1}$.\\
 
\n (ii) When $\displaystyle{ \frac{1}{2} < s - 4[\frac{s}{4}] \leq \frac{3}{2}}$, then\\
\n $\displaystyle{\phi_{k}(0) = h_{1}^{(k)}(0), \quad \phi_{k}^{'}(0) = h_{2}^{(k)}(0), \quad\quad k = 0, 1, ......, [\frac{s}{4}]-1}$:\\
\n $\displaystyle{\phi_{k}(0) = h_{1}^{(k)}(0), \quad\quad \quad\quad \quad\quad \quad\quad \quad\quad k = [\frac{s}{4}]}$.\\

\n (iii) When $\displaystyle{s - 4[\frac{s}{4}] > \frac{3}{2}}$, then\\
\n $\displaystyle{\phi_{k}(0) = h_{1}^{(k)}(0), \quad \phi_{k}^{'}(0) = h_{2}^{(k)}(0), \quad\quad k = 0, 1, ......, [\frac{s}{4}]}$.\\
\n Here,\\
$\left \{ \begin{array} {c} \displaystyle{ \phi_{0}(x) = \phi(x)},\quad\quad \quad\quad\quad\quad \quad\quad\quad\quad\quad\quad\quad\quad\quad\quad\quad\quad \quad\quad\\
\quad \quad\quad\quad\quad \quad\quad \quad\quad \\
\displaystyle{ \quad\quad \phi_{k}(x) = -\phi_{k-1}^{''''}(x) - \delta \phi_{k-1}^{'''}(x)  - \phi_{k-1}^{''}(x)  - \sum_{j=0}^{j-1} C_{k-1}^{j} \phi_{j}(x) \phi_{k-j-1}^{'}(x)} \\
\end{array}\right.$\\
\n where $\displaystyle{ C_{n}^{p} = \frac{n!}{p!(n-p)!}}$ and $[\bullet]$ is the floor function.\\
\end{definition}
\n Their main result is stated as follows.\\
\begin{theorem} ({\color{blue}[[25]]})\\
\n(i) If $s \geq 0$ and $\displaystyle{(\phi, h_{1}, h_{2}) \in H^{s}(\mathbb{R}^{+}) \times H^{\frac{s}{4} + \frac{3}{8}}(0, T) \times H^{\frac{s}{4} +\frac{1}{8} }(0, T)}$ is $s$-compatible, then equation (7.1) admits a unique solution $\displaystyle{u \in C([0, T]; H^{s}(\mathbb{R}^{+})  \bigcap L^{2}(0, T; H^{s+2}(\mathbb{R}^{+})}$ with $\displaystyle{ u_{x} \in C([0, +\infty); H^{\frac{s}{4} +\frac{1}{8} }(0, T)}$.\\
\n  Moreover, the corresponding solution map from the space of initial and boundary data to the solution space is continuous.\\
\n (ii) If $-2 < s < 0$,  $\displaystyle{(\phi \in \in H^{s}(\mathbb{R}^{+})}$ , $\displaystyle{(h_{1}, h_{2})  \times H^{\frac{s}{4} + \frac{3}{8}}(0, T) \times H^{\frac{s}{4} +\frac{1}{8} }(0, T)}$ and $\displaystyle{(t^{\frac{\vert s \vert}{4} + \epsilon}h_{1}, t^{\frac{\vert s \vert}{4} + \epsilon}h_{2})  \in H^{\frac{3}{8}}(0, T) \times H^{\frac{1}{8} }(0, T)}$  then equation (7.1) admits a unique solution $\displaystyle{u \in C([0, T]; H^{s}(\mathbb{R}^{+})}$.\\
\n  Moreover, the corresponding solution map from the space of initial and boundary data to the solution space is continuous.\\
\end{theorem}
\n {\color{red}{\bf {\large References}}}\\

\n  {\color{red}[[1]]} Adams, R.A. : Sobolev Spaces, Academic Press, New York, 1975.\\

\n    {\color{red}[[2]]}  Aimar, M.T., Etude num\'erique d'une \'equation d'\'evolution non lin\'eaire d\'ecrivant l'instabilit\'e thermodiusive d'un front de flamme (Kuramoto-Sivashinsky \'equation), PhD, University of Provence , Marseille (1982)\\

\n  {\color{red}[[3]]}  Aimar, M.T., Penel, P., Auto-organisation transitoire du mod$\grave{e}$le de Kuramoto-Sivashinsky, Lecture Notes in Physics, 241, Springer-Verlag, (1985), pp. 139-147\\

\n  {\color{red}[[4]]} Aronszajn, N. : Boundary values of functions with finite Dirichlet integral, Tech. Report of Univ. of Kansas 14 (1955) 77-94.\\

 \n{\color{red}[[5]]}  Bona, J, L and  Smith, R : The initial-value problem for the Korteweg de Vries equation. Philos. Trans. Roy. Soc London  A 278 (1975), 555-601 \\
 
 \n {\color{red}[[6]]} Cohen, B.I., Krommes, J.A., Tang, W.M., Rosenbluth, M.N.: Nonlinear saturation of the dissipative trapped-ion mode by mode coupling. Nucl. Fusion 16, 971-992 (1976)\\
 
 \n  {\color{red}[[7]]}  Demengel, F. and  Demengel, G. : Functional spaces for the theory of elliptic partial differential equations, Universitext, Springer, London ; EDP Sciences, Les Ulis, 2012.\\
 
 \n  {\color{red}[[8]]} Di Nezza, E.,  Palatucci, G, and Valdinoci, E. : Hitchhiker's guide to the fractional Sobolev spaces, Bull. Sci. Math. 136 (2012), no. 5, 521-573.\\
 
\n  {\color{red}[[9]]} Evans, L.C. : Partial Differential Equations, Graduate Studies in Mathematics, Providence, Rhode Island, American Mathematical Society, AMS, 19, 1988.\\
 
 \n  {\color{red}[[10]]} A. Friedman, A. : Partial Differential Equations, Holt, Reinhart and Winston, New York, 1969.\\

\n {\color{red}[[11]]}  Gagliardo, E. : Proprieta di alcune classi di funzioni in piu variabili, Ric. Mat. 7 (1958) 102-137.\\

\n {\color{red}[[12]]} Kato, T. : On the Korteweg-de Vries equation. Manuscr. Math. 28 (1979) 89-99.\\

\n  {\color{red}[[13]]} Kato, T. : On the Cauchy problem for the (generalized) Korteweg-de Vries equations. Adv. Math. Suppl. Stud. 8 (1983) 93-128.\\

\n {\color{red}[[14]]} Kato, T. : Quasilinear equations of evolutions with applications to partial dif- ferential equation, Lecture Notes in math.448 (1975), Springer-Verlag, 27- 50\\

\n  {\color{red}[[15]]} Kawahara, T., : Formation of saturated solitons in a nonlinear dispersive system with instability and dissipation. Phys. Rev. Lett., 51:381382, (1983).\\

\n  {\color{red}[[16]]} Kawahara,T. and Toh, S., Nonlinear dispersive periodic waves in the presence of instability and damping. Phys. Fluids, 28:16361638, (1985).\\

\n {\color{red}[[17]]} Kawahara, T. and Toh, S., Pulse interactions in an unstable dissipativedispersive nonlinear system. Phys. Fluids, 31(8):21032111, (1988)\\

\n  {\color{red}[[18]]} Kawahara, T. : Oscillatory solitary waves in dispersive media, J. Phys. Soc. Japan 33 (1972), 260-264.\\

\n  {\color{red}[[19]]} Kenig, C.,  Ponce, G. and  Vega, L. : On the (generalized) Korteweg-de Vries equation. Duke Math. J. 59 (1989) 585-610.\\

\n {\color{red}[[20]]} Kichenassamy, S and Olver, P.J. : Existence and nonexistence of solitary wave solutions to high-order model evolution equations, SIAM J. Math. Anal. 23(1992),
1141-1166.\\ 

\n {\color{red}[[21]]} Kuramoto,Y., Instability and turbulence of wave fronts in reaction diffusion systems, Prog.Theor.Phys., 63, 6, (1980) pp.1885-1903\\

\n {\color{red}[[22]]} Kuramoto, Y.: Diffusion-induced chaos in reaction systems. Prog. Theor. Phys. Suppl. 64, 346?367 (1978)\\

\n  {\color{red}[[23]]} Kuramoto, Y., Tsuzuki, T.: On the formation of dissipative structures in reaction-diffusion systems. Theor. Phys. 54, 687-699 (1975)\\

\n  {\color{red}[[24]]} Kuramoto, Y., Tsuzuki, T.: Persistent propagation of concentration waves in dissipative media far from thermal equilibrium. Prog.Theor.Phys.55, 356-369 (1976)\\

\n {\color{red}[[25]]} Li, J., Zhang, B-Y and  Zhang, Z. : A nonhomogeneous boundary value problem for the Kuramoto-Sivashinsky equation in a quarter plane, arXiv:1607.00506v2 [math.AP] 14 Jul 2016.\\

\n {\color{red}[[26]]} Manneville, P., Macroscopic Modeling of Turbulent Flows, volume 230.Lecture Notes in Physics, Springer-Verlag Berlin Heidelberg, (1985).\\

\n  {\color{red}[[27]]} Michelson, D. : Steady solutions of the Kuramoto-Sivashinsky equation, Physica D, 19 (1986), pp. 89-111.\\

\n {\color{red}[[28]]} Morgolis,S.B., Sivashinsky, G.I., Some new solutions of the flamme front equation, Lecture Notes in Physics, 241, Springer-Verlag, (1985), pp. 131-136\\

\n  {\color{red}[[29]]} Nirenberg, L. : Topics in Nonlinear Functional Analysis, Courant Institute, New York, 1974.\\

\n   {\color{red}[[30]]} Ponce, G. : Lax pairs and higher order models for water waves, J. Differential Equations 102 (1993), 360-381.\\

\n  {\color{red}[[31]]} Sell, G. R. and Taboada, M. , Local dissipativity and attractors for theKuramoto- Sivashinsky equation in thin 2D domains, Nonlin. Anal., 18, (1992), 671-687.\\

\n  {\color{red}[[32]]} Sivashinsky, G.I., Difusion-thermal theory of cellular flame, Comb.Sci.and Thech. 15, 4, (1977), pp. 137- 145\\

\n  {\color{red}[[33]]} Sivashinsky, G.I., Non linear analysis of hydrodynamic instability in laminar flames, Part I, derivation of basic equations, Acta Astronautica, 4, (1977), pp. 1117-1206\\

\n {\color{red}[[34]]} Slobodeckij, L.N. :Generalized Sobolev spaces and their applications to boundary value problems of partial differential equations, Leningrad. Gos. Ped. Inst. Ucep. Zap. 197 (1958) 54-112.\\

\n  {\color{red}[[35]]} Sobolev, S.L. : On a theorem of functional analysis, Mat. Sb. 4 (46) 1938, 39-68 (translated into English in 1963).\\

\n  {\color{red}[[36]]} Temam, R., Infinite-Dimension Dynamical Systems in Mechanics and Physics, Appl. Mat. Sc. 68, 2nd ed., Springer (1997).\\

\n  {\color{red}[[37]]} Topper, J., Kawahara, T.: Approximate equations for long nonlinear waves on a viscous fluid. J. Phys. Soc. Jpn. 44, 663-666 (1978)\\

\n {\color{red}[[38]]}Whitham, G.B. :Linear and nonlinear waves, John Wiley $\&$ Sons Inc., New York, 1999. Reprint of the 1974 original, A Wiley-Interscience Publication.\\

\n  {\color{red}[[39]]} Zhang, B.Y. ; Some  results for the nonlinear dispersive wave equation with application to control, Ph.D. Thesis, University of Wisconsin-Madison, 1990.\\
\end{document}